\documentclass[11pt,a4paper]{article}
\usepackage{amsthm,bm}
\usepackage[small]{titlesec}
\usepackage{showlabels}
\usepackage{relsize}
\usepackage{graphicx}
\usepackage[latin1]{inputenc}
\usepackage{amsmath,amssymb}
\usepackage{latexsym}
\usepackage{hyperref}
\usepackage{graphicx}
\usepackage{longtable}
\usepackage{subfigure}
\usepackage{epstopdf}
\usepackage{multirow}
\usepackage{caption}
\usepackage{hyperref}
\usepackage{fancyhdr}
\usepackage{bbm}            
\usepackage{dsfont}         
\usepackage{color}
\usepackage{float}
\usepackage{geometry}
\newtheorem{defn}{Definition}
\newtheorem{prop}{Proposition}

\usepackage[normalem]{ulem}   
\newcommand\redout{\bgroup\markoverwith
{\textcolor{red}{\rule[.5ex]{2pt}{0.4pt}}}\ULon}
\usepackage{xcolor}
\newcommand{\cub}{\textsc{cub\,}}

\begin{document}

\title{A Model-Based Fuzzy Analysis of Questionnaires \thanks{The original source of publication is available in {\sl Statistical methods and Applications}, 28:187-215 (2019) at https://link.springer.com/article/10.1007$\%2$Fs$10260$-$018$-$00443$}}

\author{Elvira Di Nardo \thanks{elvira.dinardo@unito.it, Department of Mathematics \lq\lq G. Peano\rq\rq, University of Turin, Via Carlo Alberto, 10, 10123 Torino, Italy} \and Rosaria Simone \thanks{rosaria.simone@unina.it (corresponding author) Department of Political Sciences, University of Naples Federico II, Via Leopoldo Rodin\`o, 22, 80138 Napoli, Italy}}

\date{}

\maketitle

\begin{abstract}
In dealing with veracity of data analytics, fuzzy methods are more and more relying on probabilistic and statistical techniques to underpin their applicability. Conversely, standard statistical models usually disregard to take into account the inherent fuzziness of choices and this issue is particularly worthy of note in customers' satisfaction surveys, since there are different shades of evaluations that classical statistical tools fail to catch.  Given these motivations, the paper introduces a  model-based fuzzy analysis of questionnaire with sound statistical foundation,  driven by the design of a hybrid method that sets in between fuzzy evaluation systems and statistical modelling. The proposal is advanced on the basis of \cub mixture models to account for uncertainty in ordinal data analysis and moves within the general framework of Intuitionistic Fuzzy Set theory to allow membership, non-membership, vagueness and accuracy assessments. Particular emphasis is given to defuzzification procedures that enable uncertainty measures also at an aggregated level. An application to a survey run at the University of Naples Federico II about the evaluation of Orientation Services supports the efficacy of the proposal.
\end{abstract}
\section{Introduction}
\label{intro}
Consider an experimental design aimed to investigate an unobservable trait of a population through measurements of opinions, judgements, or preferences. Then, once a questionnaire is administered and data are collected  as ratings on an ordinal scale, policy makers should not disregard to take into account the fuzziness of the outcomes. Indeed, latent constructs like customer satisfaction are inherently vague. For instance, if on a scale ranging from $1=$ `completely dissatisfied' to $7=$ `completely satisfied', the rater marks $R=6$ ($R=2$, resp.), how strong should our confidence be that he/she is actually satisfied (dissatisfied)? How confident can the scholar be about the resulting classification? In the following, we consider the satisfaction of respondents as latent phenomenon, but the proposed investigation can be applied more generally for the analysis of agreement and preferences.

A classical way to assess the imprecision and the uncertainty of evaluations is offered by Fuzzy Sets Theory \cite{Zadeh}.  When a respondent marks \lq\lq very satisfied\rq\rq about a particular service, he/she is producing a judgment on the veracity of the statement \lq\lq the quality of service is high\rq\rq,  whose plausibility can be encoded in a Fuzzy measure.  The Intuitionistic Fuzzy Set theory (IFS) \cite{Atanassov} proposes to convey the dual assessment of the membership grade  through a non-membership function. Then the residual degree of indecision and its complementary to one are considered as measures of uncertainty and accuracy, respectively, measuring the veracity of the expressed rates. Usually, fuzzy methods for evaluation are preferred over simple descriptive analysis, but they lack of a sound inferential background that could enhance its applicability. Moreover, they are strongly grounded on subjective choices, which makes reproducibility of analysis and conclusions unstable. Thus, the integration with a proper statistical tool could underpin fuzzy methods' reliability.

The present paper aims to pursue this task by fostering the application of a suitable statistical method within the fuzzy framework. Due to its psychological and probabilistic structure, an appealing candidate for our purpose is the class of \cub models \cite{DelPic2005,Ianp12}. 
This rationale conceives the data generating process driving latent perceptions into discrete evaluations as the combination of two components:  the \textit{feeling}, responsible for the level of agreement/pleasantness towards the item under investigation, and the \textit{uncertainty}, accounting for the overall nuisance affecting a fully meditated response (laziness, difficulties in understanding the question, ignorance of the topics, wording and length of the scale etc). \cub models are then defined as a two-component mixture distribution: a shifted Binomial for feeling and a discrete Uniform for uncertainty, which explains the acronym \cub: \textit{C}ombination of a \textit{Uniform} and a shifted \textit{B}inomial.

Then, the proposal is a new Fuzzy evaluation system for ratings in the setting of Intuitionistic theory in order to properly account for uncertainty in the data as meant by the \cub paradigm:  dually,   the heuristic definition of \cub uncertainty as a measure of the intrinsic fuzziness of the decision process is properly justified.  As for all fuzzy schemes \cite{Lalla}, the resulting model-based fuzzy system places itself as a support tool to the analysis of ratings.

The advantage of the proposed method is manifold: it is more objective, since fuzzy functions are structured on both data and inferential procedures and simultaneously they are designed to be more sensitive to measurement errors; it is able to discriminate among items of a questionnaire not only at an aggregated level (that is when respondents are grouped according to their choices) and it lends itself to broader extensions  to consider further  sources of the fuzziness blurring the response distribution. In particular, the resulting modeling allows us to include the so-called \textit{shelter effect}, occurring when a proportion of respondents identifies a category as a \textit{refuge} option for non-meditated choices \cite{Ian12a,IannarioPiccolo2016}, yielding to inflation of frequency in some categories. 

Our proposal considers membership and non-membership functions of spline type \cite{Marasini}, grounded on sampled information and limiting the subjectivity of parameters choice. The application of the empirical distribution function within the fuzzy evaluation system is supported by the literature \cite{Cheli,Zani} and, in these terms, it allows us to take into account also the feeling component as meant by \cub models. 

In the end, a Fuzzy analysis of questionnaire is completed with a so-called \textit{defuzzification} procedure. This last step consists in computing synthetic measures that encode a simultaneous examination of all items across respondents \cite{Leekwijck}: specifically, fuzzy functions have to be suitably weighted and aggregated to produce fuzzy composite indicators. In \cite{Marasini}, different criteria and quantification methods are discussed. In particular, the weights associated with items can be either uniform or assigned by experts who are in charge of discriminating the items by assessing their relative importance to the universe of discourse.   Here we propose an aggregator belonging to the class of Intuitionistic Weighted Aggregator Means (IWAM) \cite{Beliakov}, designed to balance  both  satisfaction and dissatisfaction for indecision  and unpredictability. 

With regard to the literature on \cub models, the novelty of the paper is the multi-item perspective offered by the defuzzification procedure, in which dimensions affected by a larger uncertainty are recognized a weaker importance. Although \cub models are designed to run an item-by-item investigation, first multidimensional perspectives are given in \cite{AndreisFerrari,Corduas}, whereas multi-object analysis are discussed in \cite{Ianp12,Capecchi} and multi-item aggregation is pursued with a model-based composite indicator in \cite{CapecchiSimone2018}. Very recently, a multivariate extension of \cub models is proposed in \cite{Colombi}.

The work is organized as follows: Section \ref{sec:secIFS} provides a short overview of the IFS theory, with a focus on composite indicators given in Section \ref{sec:FuzzyInd}. \cub models are shortly described in Section \ref{sec:CUB} and the proposed \cub-Fuzzy evaluation system is introduced in Section \ref{sec:CUBfuzzy}, where a new Fuzzy composite indicator is defined in terms of the \cub uncertainty parameter. Finally, in Section \ref{Ver}, the proposal is illustrated on the basis of a survey collected at the University of Naples Federico II about the evaluation of Orientation Services. The discussion is pursued by assuming a comparative perspective with standard methodologies for multi-item analysis. Summarizing remarks and some notes on future developments end the paper. The whole analysis has been run within the \texttt{R} environment: the code is available upon request from Authors.
\section{Fuzzy systems: Intuitionistic theory}\label{sec:secIFS}
Let $X$ be the universe of discourse. For instance, assume that we are investigating customers' satisfaction and a questionnaire is designed to that purpose: then $X$ is the set of all customers, which is observed through the respondents to the survey.  A Fuzzy set $A$ consists of a subset of $X$ endowed with a \textit{membership} function $\mu_A$ assessing the degree of membership to the set $A$,
\begin{equation*}
\mu_A :\, X \longrightarrow [0,1], \qquad x \longmapsto \mu_A(x),
\end{equation*}
in such a way that $\mu_A(x) =1$ if and only if $x$ is certainly an element of $A$, while $\mu_A(x) =0$ if and only if $x$ is certainly not.  For the illustrative example above, $A$ is the subset of the satisfied customers. Assume that answers to an item are collected on an Likert-type scale with $m=10$ categories that can be coded with equispaced integer scores, and that only categories $j\geq 6$ have a positive wording. Evaluations are not \textit{crisp}, and classification of respondents should be elastic, thus it is an over-simplification to consider certainly satisfied those users whose rate expresses satisfaction, regardless of the position along the scale. Certainly, $j=6$, $j=8$ and $j=10$ have to be associated with different degrees of belonging to $A$. A Fuzzy evaluation system will frame this circumstance by assigning increasing levels of membership to increasing scores. Then, it is widely acknowledged that in studies like those on customer satisfaction, it is of foremost importance to accompany evaluation of satisfaction with measures of unsatisfaction and dissatisfaction. In this vein, the rationale of Intuitionistic Fuzzy Sets (IFS) puts forth a theory to supply the Fuzzy analysis with a \textit{non-membership} function \cite{Atanassov}:
\begin{equation*}
\nu_A :\, X \longrightarrow [0,1], \qquad  x \longmapsto \nu_A(x),
\end{equation*}
expressing the dual assessment of the non-membership grade of an element $x$ to $A$, in such a way that if $\nu_A(x) =1$, then $x$ is certainly not an element of $A$. 

Membership and non-membership values should be defined in such a way that $0 \leq \mu_A(x) + \nu_A(x) \leq 1$ \cite{Atanassov2}, thus a measure of the residual indecision about the statement ``$x \in A$'' is given by the \textit{hesitancy degree} or \textit{Fuzzy uncertainty} function:
\begin{equation}
\label{UncFunction}
u_A(x) = 1 - \mu_A(x) - \nu_A(x).
\end{equation}

Fuzzy uncertainty function (\ref{UncFunction}) parallels the usual confidence band. Indeed the \textit{Interval-Valued Fuzzy Set} (IVFS)  for an IF singleton $ <\mu_A,\,\nu_A>$ is the function \cite{intIFS}:
$$M_A: X \longrightarrow {\mathcal B}[0,1], \qquad  x \longmapsto [\mu_A(x), 1 - \nu_A(x)]\, ,$$
where ${\mathcal B}[0,1]$ denotes the Borel set of sub-intervals of the unit interval. Then the hesitancy degree is the range of $M_A(x)=[\mu_A(x), \mu_A(x) + u_A(x)].$ For subsequent purposes, let us underline that some IFS evaluation systems first characterize the hesitancy degree and then obtain the non-membership function as:
\begin{equation}
\nu_A(x) =  1 - \mu_A(x) - u_A(x).
\label{nonmemmar}
\end{equation}

Within IFS, there are indicators aiming to summarize an item performance when aggregated among subjects, and to derive information about the latent phenomenon when aggregated among items: in particular we refer to the fuzzy score and the fuzzy accuracy \cite{Xu}. The \textit{Fuzzy score}  function \eqref{fuzzyscore} indicates how strong is the classification of membership with respect to the classification of non-membership by computing how much membership and non-membership statements are far apart, that is:
\begin{equation}
\label{fuzzyscore}
s(x) = \mu_A(x) - \nu_A(x) \in [-1,1].
\end{equation}

The \textit{Fuzzy accuracy} function, instead, measures the extent to which the fuzzy classification of membership and non-membership is encompassed: 
\begin{equation}
\label{fuzzyaccuracy}
a(x) = \mu_A(x) + \nu_A(x) = 1 - u_A(x) \in [0,1],
\end{equation}
in the sense that $a(x)=1$ denotes that no undefined state is contemplated other than membership and non-membership (conversely, $a(x)=0$ indicates a fully incomplete fuzzy statement).
\subsection{Spline fuzzy systems} \label{sec:secSpline}
When discussing the IFS framework for questionnaire analysis, a benchmark approach is the one delivered in \cite{Marasini,Marasini2,Marasini3}. Let us consider a balanced Likert-type ordinal scale with an odd number $m$ of choices, with an indifference point $i_p$ located at the middle category. This choice is convenient since the indifference point thresholds membership and non-membership grades, although the setting here established could be easily extended to scales of even length. The scale is coded into integer categories, say $1,2,\dots,m$, so that a rate $r=1$ corresponds to the most unsatisfied choice; conversely, $r=m$ corresponds to an extremely satisfied answer. \\ \indent
According to \cite{Marasini}, the spline membership function is of the type:
\begin{equation}
\label{mfspline}
\mu_A(r) =
\begin{cases}
0, & \quad 1\leq r < a,\\
\dfrac{1}{2} - \dfrac{1}{2}\bigg(2\dfrac{i_p - r}{b - a}\bigg)^{\epsilon},  &  \quad  a  \leq  r \leq i_p, \\
\dfrac{1}{2} + \dfrac{1}{2}\bigg(2\dfrac{r - i_p}{b - a}\bigg)^{\epsilon},  &  \quad  i_p \leq  r \leq b, \\
1, &\quad b < r \leq m,
\end{cases}
\end{equation}
with $\epsilon >0$, and $b - a$ denoting the range of non-crisp responses (notice that for the sake of the subsequent discussion, we define membership directly on the ordinal scale, whereas the approach in \cite{Marasini} works on the latent continuous measurement scale). Spline parameters $\epsilon, \eta, \theta$ are chosen according to the sampling experiment, the strength and vagueness of the wording of the scale and its length. For instance, in \cite{Marasini} the authors advocate to adopt a linear spline ($\epsilon = 1$) for items within a certain section of the questionnaire, and a quadratic spline ($\epsilon = 2$) for items within another one, measured on a scale whose wording can be perceived vaguer in the central part of the scale, or in cases in which there is a non-linear step between subsequent categories. In the latter cases, at least a quadratic spline should be recommended.\\
\indent
The hesitancy degree is defined from \eqref{mfspline} with:
\begin{equation}
\label{uncMarasini}
u_A(r) = \mu_A(r)^{\theta}( 1-\mu_A(r))^{\eta}, \quad \theta, \eta \geq 1,
\end{equation}
and then the non-membership function is derived according to \eqref{nonmemmar}. This definition is meant to convey both membership and its residual assessment to uncertainty measurement: for balanced scale, one sets $\theta = \eta$. Then, non-membership degrees are obtained from \eqref{nonmemmar}. As a result,  without any prior experts' assessments on the values of parameters, the Fuzzy spline functions \eqref{mfspline} and \eqref{nonmemmar} are equal for all items evaluated on a common scale.

 Although interesting, we will not rely on Definition \eqref{uncMarasini} for the uncertainty function, but implement a fuzzy uncertainty that carries a specific statistical interpretation in the spirit of freeing the stakeholders from preliminary subjective assessments. 
\subsection{Fuzzy composite indicators} \label{sec:FuzzyInd}
Consider a latent phenomenon to be measured by $K$ observable variables, as the items of a questionnaire. Assume the questionnaire has
been filled out by $n$ respondents, who have chosen among $m$ ordered alternatives, available for each item. Let $\mathbf{r}_j =(r_{j,1},r_{j,2},\dots, r_{j,K})$ be the row vector of ratings given by the $j$-th respondent for $j=1,2,\ldots,n$, to the $K$ items, and denote with $\mu^{(k)}_A(\cdot), \nu^{(k)}_A(\cdot)$ the membership and non-membership functions for the $k$-th item, respectively. If seeking for a composite fuzzy value for each respondent, the IWAM (Intuitionistic Weighted Aggregator Mean) is defined as the pair:
\begin{equation}
\label{IWAM}
<\mu_A(\mathbf{r}_j), \nu_A(\mathbf{r}_j) >  \;= \; <  \sum_{k=1}^K w_k\, \mu^{(k)}_A(r_{j,k}), \,\sum_{k=1}^K w_k \,\nu^{(k)}_A(r_{j,k}) >,
\end{equation}
where $\{w_1, \ldots, w_K\}$ is a given system of weights such that $\sum_{k=1}^K w_k=1$, establishing the relative importance of items. Such values could be used to perform a fuzzy clustering of responses, where belonging of each observation to a cluster is decided on the basis of the fuzzy composite score $\mu_A(\mathbf{r}_j) - \nu_A(\mathbf{r}_j)$, for instance.
Each aggregated value is considered an IFS singleton $ <j, \mu_A(\mathbf{r}_j), \nu_A(\mathbf{r}_j) >$, thus a final composite score can be obtained by considering uniform weights for subjects:
\begin{equation}
\label{score}
<\bar{\mu},\bar{\nu}> \;=\; <\frac{1}{n}\sum_{j=1}^n \mu_A(\mathbf{r}_j),  \frac{1}{n}\sum_{j=1}^n \nu_A(\mathbf{r}_j) >.
\end{equation}
Then, according to \eqref{UncFunction} and \eqref{score}, the \textit{uncertainty} (or \textit{hesitancy degree}) is computed as the global residual degree of indeterminacy of the fuzzy assessment:
\begin{equation}
\label{uncertaintyScore}
\bar{u} = 1 - \bar{\mu} -\bar{\nu},
\end{equation}
whereas the overall Fuzzy score and Fuzzy accuracy are given respectively by:
$$
\bar{s} = \bar{\mu} - \bar{\nu}, \qquad \bar{a} = \bar{\mu}     + \bar{\nu}.
$$

Different choices of weights in \eqref{IWAM} give different indicators. In the framework of composite indicators, the choice of a weighting system is of primary importance. As a matter of fact, several applications suggest to choose weights depending on the loadings of the first principal component or factor. Nevertheless, such choice is consistent only if that variable explains a large proportion of the variability.  For this reason, and aiming to a fuzzy system that is model-based and thus not subjective, we will propose a system of weights that is driven by data through estimation procedure: in this sense, it can be considered a safer option.

\section{CUB models}\label{sec:CUB}
Let $R$ be the rating random variable modelling the response distribution to an item of a questionnaire, measured on a scale with $m$ ordered categories coded as integers from $1$ up to $m$. A \cub distribution $\cub(\pi,\xi)$ for $R$ consists in the following two-component mixture with parameters $(\pi, \,\xi) \in (0,1] \times [0,1]$:
$$
{\mathbb P}\big(R=r \mid \pi, \xi \big) =\pi\,b_r(\xi)+(1-\pi)\,h_r\,, \quad r=1,2,\dots,m\,,
$$
where $b_r(\xi)$, $r=1,2,\dots,m$ for $m>3$ denotes the shifted Binomial distribution with parameter $1-\xi$:
$$
b_r(\xi) = \binom{m-1}{r-1}\xi^{m-r}(1-\xi)^{r-1}, \quad r=1,2,\dots,m\,,
$$
and $h_r = \dfrac{1}{m}$ is the discrete Uniform distribution over the given support. The parameter $\xi$ is referred to as the \textit{feeling} parameter since $1-\xi$ measures the preference of a category over the lower ones in a sequence of pairwise comparisons among categories  \cite{Del00a,IannarioPiccolo2016}. The parameter $\pi$, instead, is called the \textit{uncertainty} parameter since $1 - \pi$ charges for the inherent fuzziness arising when perception translates into an evaluation, and thus measures the overall uncertainty of the respondent's assessment. The role of the uncertainty component within the \cub rationale has been usually considered as an expression of the inherent indeterminacy of human decisions, generating fuzziness and thus representing a source of unpredictability of the evaluation process. As a by-product and since the Uniform distribution represents the least informative model, its weight in the mixture aims at catching the level of heterogeneity in the data.\\ \indent

It should be emphasized that the choice of the Uniform distribution for the uncertainty component adheres to the baseline \cub paradigm: under this assumption, $\pi$ is an inverse indicator of heterogeneity. Departing from the defining specification, other choices can be supported to model uncertainty in the data:  response styles and category-specific measurement errors can be suitably specified by adjusting this component (see \cite{Gottard2014,SimoneTutz}). These extensions do not affect the distinctive trait of the \cub fuzzy evaluation system since this is grounded on the mixing weight $\pi$ for the deliberate choice. Alternative distributions for the uncertainty component would simply change interpretation of results and penalize data for response styles or more specific form of uncertainty. In the following, we will focus on some particular circumstances for illustrative purposes. 

For instance, in order to further disentangle the fuzziness charged by the uncertainty component, one may contemplate a \textit{shelter effect} concentrated at category $c \in \{1,\dots,m\}$ in the \cub mixture distribution when inflation in $c$ is observed \cite{Ian12a}. Let us consider a degenerate random variable $D_r^{(c)}$ such that ${\mathbb P}\bigl(D_r^{(c)}=r\bigr) = 1$ if $r=c$ and ${\mathbb P} \bigl(D_r^{(c)}=r\bigr)=0$ otherwise. Then, the \cub distribution $\cub(\pi,\xi,\delta)$ with \textit{shelter effect} at $r=c$ is:
\begin{equation}\label{shelter2}
{\mathbb P}(R=r \mid \pi^{\star}, \xi, \delta) = \delta\,D_{r}^{(c)}\,+\,(1-\delta)\,\bigl[\,\pi^{\star}\,b_{r}(\xi)+(1-\pi^{\star})\,\, h_r \,\bigr], \,\,\, r=1,2, \dots,m,
\end{equation}
for $m>4$. The additional parameter $\delta$ quantifies the importance of the \textit{shelter effect}. 
When testing its significance, it may be useful to deal with \eqref{shelter2} according to the following equivalent parameterization:
$$
{\mathbb P}(R=r \mid \pi_1, \pi_2, \xi) =  \pi_ 1 \,b_{r}(\xi) + \pi_2 h_r + (1- \pi_1 - \pi_2)\,D_{r}^{(c)}, \quad r=1,2, \dots,m,
$$
with $\pi_1 = \pi^{\star}(1-\delta) >0, \pi_2=(1-\pi^{\star})(1-\delta) \geq 0$. When the inclusion of a \textit{shelter effect} in the model yields to a significant improvement of the fit (to be checked with a Likelihood Ratio Test, for instance), the overall level of inaccuracy has to convey both the heterogeneity accounted by the Uniform distribution and the \textit{shelter effect} as measured by the parameter $\delta$. Since $\pi_1$ is the mixture coefficient corresponding to a deliberate choice, the measure of the overall uncertainty in this augmented case corresponds to $1-\pi_1 = \pi_2 + \delta$ taking into account the \textit{shelter effect}. In order to provide a general framework not limited to cases where the shelter is significant, we shall use the notation $\pi$ in place of $\pi_1$, since in that cases the whole discussion holds for baseline \cub models.
Notice that \cub models paradigm assumes a linear step between adjacent categories: for non-linear versions, see \cite{nlcub}.

For our computation we shall rely on the Maximum Likelihood (ML) estimates $\hat{\pi},\hat{\xi},\hat{\delta}$ of $\pi,\xi,\delta$, respectively (equivalently $\hat{\pi}_1, \hat{\pi}_2, \hat{\xi}$) obtained by running the Expectation-Maximization algorithm \cite{DelPic2005,Piccolo2006} as implemented in the \texttt{R} package \cub 
 \cite{CUBpackage}.
\subsection{CUB-Fuzzy evaluation system}\label{sec:CUBfuzzy}
The idea to use \cub model parameters in computing membership functions stems from the preliminary work  \cite{FuzzyCUB}, but the method here presented is more accurately designed.

For a preliminary investigation and comparison with the methods introduced in Section \ref{sec:secIFS}, we have focussed on balanced Likert-type scales of odd length with indifferent point at the midpoint.
 Suppose the scale is oriented in such a way that \lq\lq the greater the score, the higher the feeling\rq\rq, that is, there is a positive relation between the latent phenomenon and the scale. Here, {\it negative} and {\it positive} refer to expression of satisfaction, so that $r< i_p$ ($r > i_p$) corresponds to a {\it negative} ({\it positive}) evaluation.

\begin{defn}  \label{CUBspline}
For a given item of the questionnaire, the \cub-Fuzzy membership function is:
$$
\label{CUBsplineeq}
\mu_A(r) =
\begin{cases}
0, & \quad 1 \leq  r \leq l_b,\\
\dfrac{\hat{\pi}}{2}  -  \dfrac{\hat{\pi}}{2}\,\dfrac{F(i_p) - F(r)}{F(i_p) - F(l_b)}, &  \quad  l_b + 1 \leq  r \leq i_p, \\
\dfrac{\hat{\pi}}{2}  +  \dfrac{\hat{\pi}}{2}\,\dfrac{F(r) - F(i_p)}{F(u_b-1) - F(i_p)}, &  \quad  i_p \leq  r \leq u_b - 1, \\
1, & \quad u_b \leq r \leq m,
\end{cases}
$$
where $F(r)$ denotes the empirical distribution function of the given variable, $\hat{\pi}$ is estimated from a \cub model fitted to the data\footnote{$\hat{\pi}$ is replaced by $\hat{\pi}_1$ when \textit{shelter effect} is considered.} and $l_b$ ($u_b$, resp.) is a fixed lower (upper, resp.) bound to threshold the categories corresponding to crisp negative (positive, resp.) scores.\\
\end{defn}

In full generality, the setting of $l_b$ and $u_b$ may be affected by the wording of the scale, the problem under investigation and/or a preliminary analysis of the data. For $l_b = 1$ and $u_b = m$, the membership function \eqref{Zani} corresponds to the totally fuzzy and relative approach given in \cite{Cheli}.  This choice allows us to penalize uniformly each category and it is best-suited for our purpose of accounting for heterogeneity, and thus it is the natural choice for the \cub-Fuzzy proposal.\\
 
Definition  \ref{CUBspline} is a linear spline in the distribution function.  Specifically and compared with \eqref{mfspline}, definition \ref{CUBspline} relies on the \cub uncertainty parameter, but also considers a spline transformation of the empirical distribution function for the item rather than of ordinal categories as in \eqref{mfspline}. Indeed, measuring distances between categories via their differences may be inappropriate since results depend on the chosen scores: this issue is particularly relevant when the same latent trait is assessed in different groups, locations or times for comparison purposes. 
Most importantly, it is not necessary to specify spline degrees $\epsilon$, since $\hat{\pi}$ will charge for all the unspecified effects and vagueness of the evaluation, as that derived from the nature of the scale \cite{IannarioEJASA}.\\

The rationale and probabilistic genesis behind Definition \ref{CUBspline} can be summarized as follows:
\begin{description}
\item[{\it i)}] the updating of the category $r$ is penalized with the mixing weight for the feeling component $\hat{\pi}$ since it establishes the accuracy of the preference part of the model by adjusting its importance for heterogeneity and diverse sources of imprecision in the data;
\item[{\it ii)}] for $r> i_p$ ($r < i_p$, resp.) the frequency of the category $r$ is normalized taking into account the set of positive (negative, resp.) non-crisp choices;
\item[{\it iii)}] the greater is the heterogeneity (that is, as $\hat{\pi} \rightarrow 0$), the less meaningful is the contribution of the relative frequencies to the membership degrees.
\end{description}

The choice of normalizing the updating contribution $F(r) - F(i_p)$ with $F(u_b-1) - F(i_p)$ for the categories $i_p \leq r \leq u_b-1$ can be explained as follows: the categories $r \geq u_b$ are certainly associated with membership to $A$ (in our case, $A$ is the set of satisfied users) as $\mu_A(r)=1$. Symmetric arguments apply to the choice $F(i_p) - F(l_b)$ for lower categories $l_b + 1 \leq \,r \leq i_p$. Hence, the shades of membership across intermediate positive categories should rather be computed starting from the indifference point and excluding the categories being assigned crisp membership degrees. Moreover, the choice of halving $\hat{\pi}$ and distinguishing between left and right non-crisp sides of the scale is due to weight for the dual contribution of each category to the assessment of membership and non-membership. \\

The \cub-Fuzzy proposal stems from the central idea of giving to $1-\hat{\pi}$ a proper definition as measure of fuzziness of the decision process. Thus we assume the range of the IVFS constantly equal to $1-\hat{\pi}$ for each category $r$, in agreement with the role that uncertainty plays in \cub models.
\begin{defn} \label{unc}
For a given item of the questionnaire, the \cub-Fuzzy uncertainty function for $A$ is defined as:
$$
u_A(r) = \left\{ \begin{array}{ll}
0, & \quad 1 \leq r \leq l_b \,\, \hbox{and} \, \,\, u_b \leq r \leq m, \\
1 - \hat{\pi}, & \quad l_b + 1 \,\leq\, r \leq \,u_b-1.
\end{array} \right.
$$
\end{defn}
From \eqref{fuzzyaccuracy}, the Fuzzy accuracy function results to be:
$$
a(r) = \left\{ \begin{array}{ll}
1, &  \quad 1 \leq r \leq l_b \,\, \hbox{and} \, \,\, u_b \leq r \leq m, \\
\hat{\pi}, &  \quad l_b + 1 \,\leq\, r \leq \,u_b-1,
\end{array} \right.
$$
catching the propensity to assume a meditated response mechanism. Indeed, given the mixture definition, $\pi$ is a direct indicator of reliability of predictions under the feeling component, which could be adjusted to incorporate also overdispersion \cite{Iannario2014} or a more general specification \cite{Tutzb}. Thus, the choice for $u(r) = 1-\pi$ under the \cub-Fuzzy system implies that the assessment of membership and non-membership of score $r$ is penalized by the unpredictability of responses under the feeling model. In addition, $\pi$ can be interpreted as a measure of propensity between a well-structured response-behaviour and a random choice: the close $\pi \rightarrow 1$, the stronger the frequency distribution can be legitimately used for a fuzzy evaluation system.\\

From \eqref{nonmemmar}, the non-membership function $\nu_A(r)$  is given by:
\begin{equation}
\label{nmfspline}
\nu_A(r) =
\begin{cases}
1, & \quad  1 \leq r \leq l_b, \\
\dfrac{\hat{\pi}}{2} + \dfrac{\hat{\pi}}{2} \,\dfrac{F(i_p) - F(r)}{F(i_p)-F(l_b)},   & \quad l_b + 1 \leq \,r \leq i_p,  \\
\dfrac{\hat{\pi}}{2} - \dfrac{\hat{\pi}}{2} \,\dfrac{F(r) - F(i_p)}{F(u_b-1)-F(i_p)},&  \quad  i_p < r \leq u_b-1, \\
0,  &\quad u_b \leq r \leq m.
\end{cases}
\end{equation}
Let us remark that, as the \cub-Fuzzy uncertainty decreases (that is, the more $\hat{\pi}$ approaches $1$), the more the non-membership function
\eqref{nmfspline} increases towards $1$ by moving from the indifference point to the first category, and similarly decreasing in the opposite direction of the scale.  If the scale orientation is opposite, then the  definition of membership and non-membership should be simply switched.

The middle point of the scale is then equally mirrored both in the membership and non-membership scores, as
$$
\mu_A(i_p) = \nu_A(i_p)= \dfrac{\hat{\pi}}{2}.
$$
In this way, the indifference expressed by the respondent choosing $i_p$ corresponds to equi-preference of categories, since $\pi$ is an inverse indicator of heterogeneity. Then, for each rating, the degrees of membership and non-membership are equally split around the indifference point $i_p$, by halving the weight of the uncertainty parameter. Note that for distributions with low heterogeneity and thus with higher concentration, one has that $\pi \rightarrow 1$ and $\mu_A(i_p) = \nu_A(i_p) \rightarrow 1/2$, as for the spline approach recalled in Section \ref{sec:secSpline}. In view of the defuzzification procedure, the membership and non-membership degrees are defined in such a way that the accuracy is lower for the items affected by higher heterogeneity, regardless of the level of feeling. Indeed, for increasing heterogeneity (that is, as $\hat{\pi} \rightarrow 0$), from Definition \ref{CUBsplineeq} and \eqref{nmfspline} we have:
$$\mu_A(r), \nu_A(r) \rightarrow 0, \,\, u_A(r) \rightarrow 1, \,\, a(r) \rightarrow 0, \qquad r=l_b +1,\dots, u_{b}-1\, , $$
so that the residual fuzziness $u_A(r)$ increases over the accuracy $a(r)$; accordingly, we are let to negligible membership/non-membership values.\\

The usage of the empirical distribution function for a Fuzzy evaluation system is also the key of the approach pursued in \cite{Zani}, which accomplishes a questionnaire analysis in a standard Fuzzy Sets (FS) framework, grounded solely on the membership function:
\begin{equation}
\label{Zani}
\mu_A(r) =
\begin{cases}
0, & \quad 1 \leq r \leq l_b, \\
\mu_{A}(r-1) + \dfrac{F(r) - F(r-1)}{1-F(l_b)}, &  \quad  l_b < r < u_b, \\
1, &\quad u_b \leq r \leq m.
\end{cases}
\end{equation}
In the forthcoming discussion, this classical FS method will be referred to as the \textit{empirical} Fuzzy system. \\

\subsection{Scoring uncertainty}\label{sec:uncertainty}
For the \cub-Fuzzy evaluation system, we propose the IWAM \eqref{IWAM} as aggregation index, but with weights $\{w_k\}$ depending on the \cub uncertainty parameter. This choice meets the well-acknowledged recommendation to assign weights that are larger for the more explanatory items, as in \cite{Marasini}.  Here, explanatory is meant as related to accuracy in the assessment of the fuzzy trait and it is inversely related to uncertainty. Thus,  the rationale of the \cub-Fuzzy evaluation system is to penalize items with higher estimated heterogeneity, them being 
less reliable and explanatory for the assessment of membership and non-membership to $A$.  In this regard, we shall consider the Fuzzy proportion of uncertainty function \eqref{UncFunction}:
\begin{equation}
\label{FuzzyPropUF}
g(X_k) = \frac{1}{n}\sum_{j=1}^n u_A^{(k)}(r_{j,k}),\,\, \hbox{for}\,\,\,\, k=1, \ldots, K.
\end{equation}
and apply an inverse transform to impute low weight to more uncertain items:
\begin{equation}
\label{pesiZani}
w_k = \ln\bigg( \frac{1}{g(X_k)}\bigg) \bigg/ \sum_{l=1}^{K} \ln\bigg( \frac{1}{g(X_l)}\bigg), \,\, \hbox{for}\,\,\,\, k=1, \ldots, K.
\end{equation}
Here the logarithm transform is taken only to prevent excessive values for  very low uncertainty.  This weighting scheme has been already used in the Fuzzy Sets literature (that is, only with reference to membership functions) \cite{Zani,Zani2} to assess the capabilities of each category $r$ in expressing satisfaction across items:
\begin{equation}\label{(crispZani)}
\tilde{\mu}_A(r) = \sum_{k=1}^{K} w_k\, \mu^{(k)}_A(r)\, , \quad r=1,2,\dots,m.
\end{equation} 
In that case, the weights have been based on the Fuzzy proportion of the achievement of the target (in our case, respondents' satisfaction):
\begin{equation}
\label{FuzzyProp}
g(X_k) = \frac{1}{n}\sum_{j=1}^n \mu^{(k)}_A(r_{j,k}), \,\, \hbox{for}\,\,\,\, k=1, \ldots, K.
\end{equation}
for which the transformation \eqref{pesiZani} prevents from giving higher importance to the \textit{rare features} among subjects.


More generally, as we are considering that all items are collected on the same ordinal scale, from Definition \ref{unc} with $r$ replaced by $r_{j,k},$ the proportion $g(X_k)$ in \eqref{FuzzyPropUF} has the following closed form.

\begin{prop}
If $\hat{\pi}^{(k)}$ is the estimated \cub uncertainty parameter of the $k$-th item, then
\begin{equation}
\label{FuzzyPropUF1}
g(X_k) =  \big(1 - \hat{\pi}^{(k)}\big) \left( F^{(k)} \big(u_b^{(k)}-1 \big) - F^{(k)} \big(l_b^{(k)}\big)\right),\,\, \hbox{for}\,\,\,\, k=1, \ldots, K,
\end{equation}
where $F^{(k)}(\cdot)$ is the empirical distribution function of ratings on the $k$-th item.\\
\end{prop}

Note that $F^{(k)} \big(u_b^{(k)}-1 \big) - F^{(k)} \big(l_b^{(k)}\big)$ is the percentage of respondents for which $l_b^{(k)} < r_{j,k} < u_b^{(k)}$, thus whose fuzzy evaluation on the $k$-th item is not crisp. 

If $l_b^{(k)}=1$ and $u_b^{(k)}=m,$ equation (\ref{FuzzyPropUF1}) simplifies in the \cub-Fuzzy uncertainty function $\bar{u}^{(k)}$:
\begin{equation}
\bar{u}^{(k)} = \dfrac{1}{n}\sum_{j=1}^n u^{(k)}_A(r_{j,k}) =  \big(1 - \hat{\pi}^{(k)}\big) \left( F^{(k)}(m-1) - F^{(k)}(1) \right) 
\label{unck}
\end{equation}
 when aggregating the $k$-th item among respondents. Then
the Fuzzy uncertainty score $\bar{u}$ in \eqref{uncertaintyScore} can be written as:
$$
\bar{u} = \sum_{k=1}^K w_k \bar{u}^{(k)} = \sum_{k=1}^{K} w_k \big(1 - \hat{\pi}^{(k)}\big) \left( F^{(k)} \big(m-1 \big) - F^{(k)} \big(1\big)\right).
$$
In this sense, the \cub uncertainty parameters are given a precise fuzzy interpretation also at the aggregated level. 

\section{A case study}
\label{Ver}
Fuzzy methods for questionnaire analysis are particularly appealing in evaluation studies. Motivated by this feature, we show how the \cub-fuzzy proposal can be applied on the assessment of satisfaction for the Orientation Services at University of Naples Federico II. 
The survey was administered from 2002 to 2008 across all the 13 Faculties  and aimed at measuring the satisfaction towards the service\footnote{Data are available at \text{\url{http://www.labstat.it/home/research/resources/cub-data-sets-2/}.}} across different dimensions of the trait. On a balanced $m=7$ point Likert scale: $1$ = `extremely unsatisfied', $2$ = `very unsatisfied', $3$ = `unsatisfied', $4$ = `indifferent', $5$ = `satisfied', $6$ = `very satisfied', $7$ = `extremely satisfied', the following measurements were collected:
\begin{itemize}
\item satisfaction on the acquired information (\texttt{informat});
\item evaluation of the willingness of the staff (\texttt{willingn});
\item adequacy of time-table of opening-hours (\texttt{officeho});
\item evaluation of the competence of the staff (\texttt{compete});
\item global satisfaction (\texttt{global}).
\end{itemize}
The present discussion will concern the data collected in 2002, consisting of $n=2179$ observations. Motivations for our choice include the fact that the evaluation of University courses, offices and institutions is a popular topic for fuzzy analysis; in addition, the first available wave was chosen since these data allows us to discuss and illustrate all the nuances of the proposal.

In the first part of the section, \cub models  are fitted to the data: the estimation procedure takes into account the \textit{shelter effect}, if significant. Then, the \cub-Fuzzy system introduced in Section \ref{sec:CUBfuzzy} is compared with the empirical and the spline approaches recalled in Section \ref{sec:CUBfuzzy} and \ref{sec:secIFS} respectively, within the classical and IF settings, respectively.
\subsection{CUB models estimation}
The ML estimates of \cub parameters for the chosen data are summarized in Table \ref{tab:param} \footnote{Estimation and test for \cub models are run via the \texttt{R} package `CUB' available on CRAN.}. Overall, there is a moderate level of uncertainty and an extreme positive feeling across the items, the highest satisfaction being expressed for \texttt{willingn} ($1-\hat{\xi}= 0.8833$), the lowest for \texttt{officeho} ($1-\hat{\xi}= 0.8029$). However, there are certain items for which uncertainty is not negligible and an evaluation system should properly consider these differences. In particular, \texttt{officeho} is the item with the highest estimated uncertainty ($1-\hat{\pi}=0.3198$), followed by \texttt{informat} and then by \texttt{compete}. \\
\begin{table}[h!]
\caption{Parameter estimates: \cub model no shelter effect (standard errors in parentheses)}
\label{tab:param}     
\centering
\scalebox{1}{\begin{tabular}[pos=c]{p{2cm}p{1.5cm}p{1.5cm}p{1.5cm}p{1.5cm}p{1.5cm}}
\hline \hline 
&  \texttt{informat} & \texttt{willingn} & \texttt{officeho} & \texttt{compete} &  \texttt{global} \\
\hline  $\hat{\pi}$  & $\underset{\phantom{-}\emph{(0.016)}}{0.794}$ & $\underset{\phantom{-}\emph{(0.012)}}{0.857}$ & $ \underset{\phantom{-}\emph{(0.019)}}{0.680} $&$\underset{\phantom{-}\emph{(0.015)}}{0.802} $& $\underset{\phantom{-}\emph{(0.013)}}{0.868} $\\
$\hat{\xi}$  &   $\underset{\phantom{-}\emph{(0.005)}}{0.181}$ & $\underset{\phantom{-}\emph{(0.004)}}{0.117} $& $\underset{\phantom{-}\emph{(0.006)}}{0.197}$ & $\underset{\phantom{-}\emph{(0.005)}}{0.164} $& $\underset{\phantom{-}\emph{(0.004)}}{0.171}$  \\
\noalign{\smallskip}\hline \noalign{\smallskip}
\end{tabular}}
\end{table} 

Table \ref{tab:param_she} shows the estimation results when including shelter effects: globally, the previous comments continue to hold, but uncertainty is disclosed in more details and feeling estimates are corrected. 
\begin{table}[h!]
\caption{Parameter estimates: \cub model with shelter at category $c$ (standard errors in parentheses)}
\label{tab:param_she}       
\centering
\scalebox{0.9}{\begin{tabular}[pos=c]{p{3.2cm}p{1.5cm}p{1.5cm}p{1.5cm}p{1.5cm}p{1.5cm}}
\hline \hline 
&  \texttt{informat} ($c=5$)    & \texttt{willingn} ($c=7$) & \texttt{officeho} ($c=7$) & \texttt{compete} ($c=7$)&  \texttt{global} ($c=7$) \\
\hline \hline $\hat{\pi}^{\star}$      & $\underset{\phantom{-}\emph{(0.017)}}{0.749}$ &$\underset{\phantom{-}\emph{(0.012)}}{0.848}$ & $\underset{\phantom{-}\emph{(0.022)}}{0.705}$& $\underset{\phantom{-}\emph{(0.016)}}{0.804 }$& $\underset{\phantom{-}\emph{(0.013)}}{0.869 }$\\
$\hat{\xi}$                     &   $\underset{\phantom{-}\emph{(0.006)}}{0.153}$& $\underset{\phantom{-}\emph{(0.005)}}{0.157}$ &$\underset{\phantom{-}\emph{(0.009)}}{0.265 }$ & $\underset{\phantom{-}\emph{(0.008)}}{0.199 }$ & $\underset{\phantom{-}\emph{(0.007)}}{0.183}$ \\  
$\hat{\delta}$                  &   $\underset{\phantom{-}\emph{(0.014)}}{0.085}$ &$\underset{\phantom{-}\emph{(0.011)}}{0.194}$ & $\underset{\phantom{-}\emph{(0.014)}}{0.145}$ & $\underset{\phantom{-}\emph{(0.020)}}{0.117}$ & $\underset{\phantom{-}\emph{(0.022)}}{0.048}$\\
\hline \hline $\hat{\pi}_1$ (Accuracy) & $\underset{\phantom{-}\emph{(0.021)}}{0.685 }$& $\underset{\phantom{-}\emph{(0.015)}}{0.683}$  &$\underset{\phantom{-}\emph{(0.019)}}{0.602 }$ &$\underset{\phantom{-}\emph{(0.021)}}{0.710 }$ &$\underset{\phantom{-}\emph{(0.022)}}{0.828}$ \\ 
$\hat{\pi}_2$                   &   $\underset{\phantom{-}\emph{(0.015)}}{0.230}$ & $\underset{\phantom{-}\emph{(0.012)}}{0.123}$  &$\underset{\phantom{-}\emph{(0.019)}}{0.252 }$ &$\underset{\phantom{-}\emph{(0.015)}}{0.173}$  &$\underset{\phantom{-}\emph{(0.013)}}{0.125}$\\ \hline \hline
$ 1-\hat{\pi}_1$ (Uncertainty)  &  0.315 & 0.317 &0.398 & 0.290 & 0.172 \\
$ 1-\hat{\xi}$ (Feeling)  & 0.847 & 0.843 &0.735 &0.802 &0.817 \\
\hline \hline 
\end{tabular}}
\end{table}
In particular, notice that \texttt{informat} has a different shelter category ($c=5$) compared with all other items (for which $c=7$).   Despite the shelter category is the same for the last four items, its effect is more prominent for \texttt{willingn} ($\hat{\delta} =0.194$), associated also with the lowest weight of the Uniform distribution ($\hat{\pi}_2= 0.123$). Instead, even if significant, the \textit{shelter effect} within \texttt{global} is the weakest; moreover, this item corresponds to the strongest attitude towards a more meditated choice, resulting in the highest level of accuracy ($\hat{\pi}_1 = 0.828$). For the sake of brevity, statistical results on the model selection (based on the BIC criterion) are skipped and are available on demand. \\

The rationale of the \cub-Fuzzy evaluation system is to provide statistical models for rating data with some veracity analytics in the spirit of a fuzzy analysis of questionnaire. Then, the proposal is not advanced to be a direct competitor of ordinary techniques; nevertheless, it is worth to notice how its performances match with outcomes of standard procedures. 

In this respect, the first two components as obtained from a PCA are sufficient for our purposes since they account for more than 80\% of the total variability. 
Results confirm that item \texttt{officeho} plays a distinctive role in the assessment of the latent satisfaction, thus it should be properly discriminated and weighted.
\begin{figure}[h!]
\centering
\includegraphics[scale=0.4]{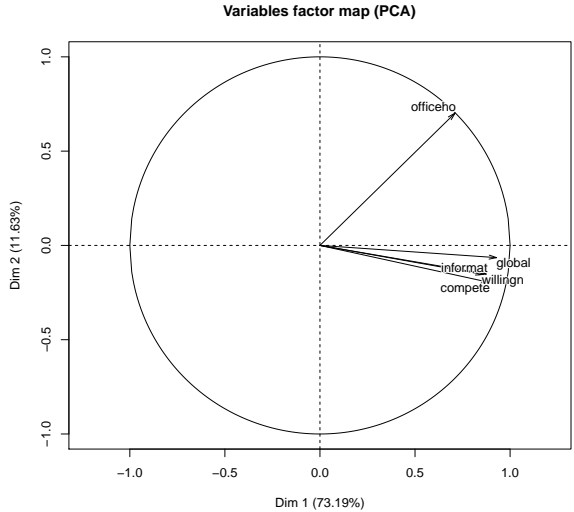}
\caption{PCA: Variable factor map.}
\end{figure}

Notice that Likert-scale categories are ordered measurements that can be thought of as cutpoints of a latent continuum. Thus, it is not always adequate to fit a PCA or Factor analysis directly on a data matrix like that from a rating survey, unless a proper correlation is obtained. In this regard, polychoric correlation is a validated choice especially when the number of categories is moderate. 

\subsection{Classical FS: CUB-Fuzzy versus empirical}\label{sec:CUB_Zani}
\begin{figure}[h!]
\centering{\includegraphics[scale=0.3]{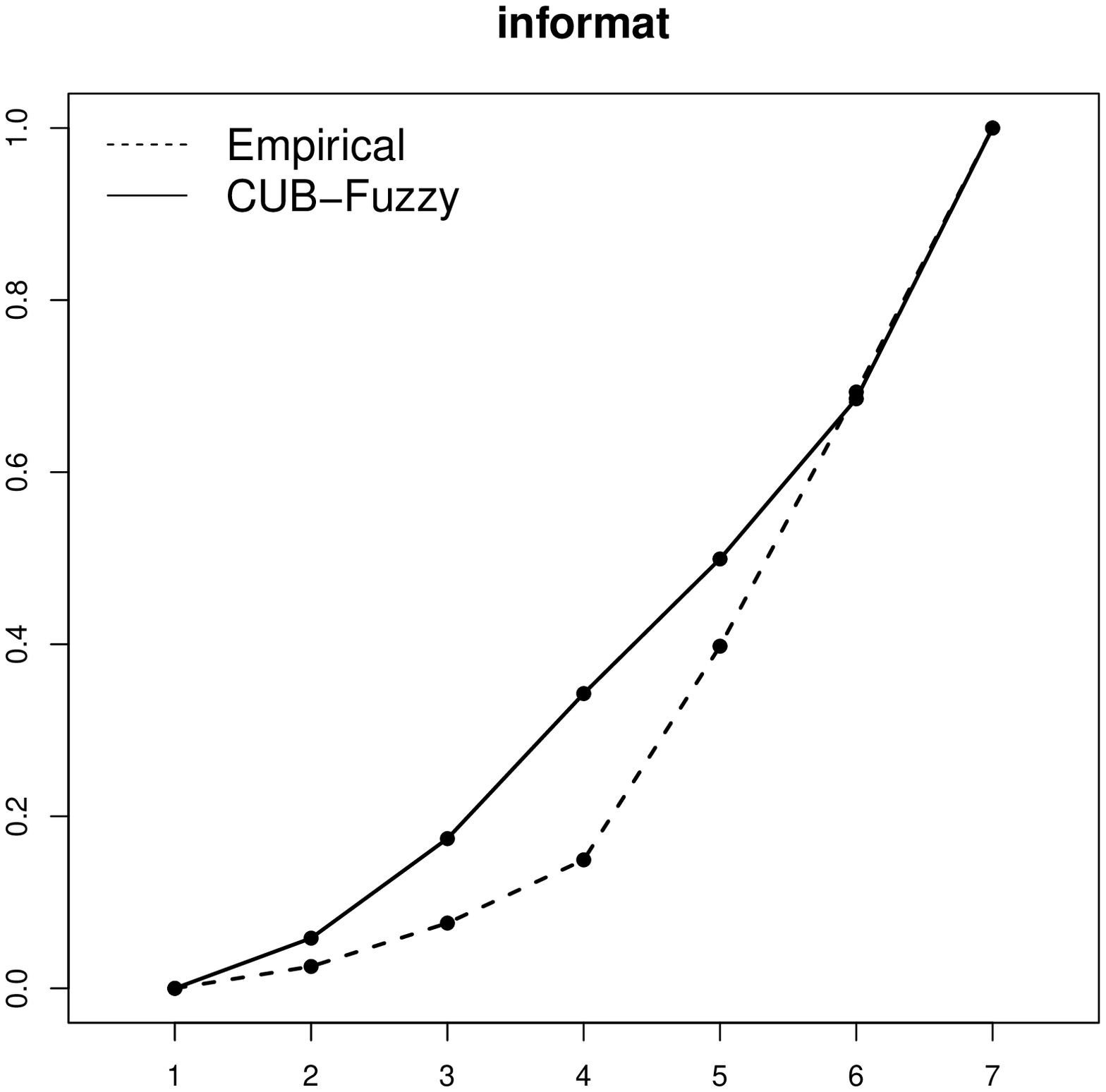}
\includegraphics[scale=0.30]{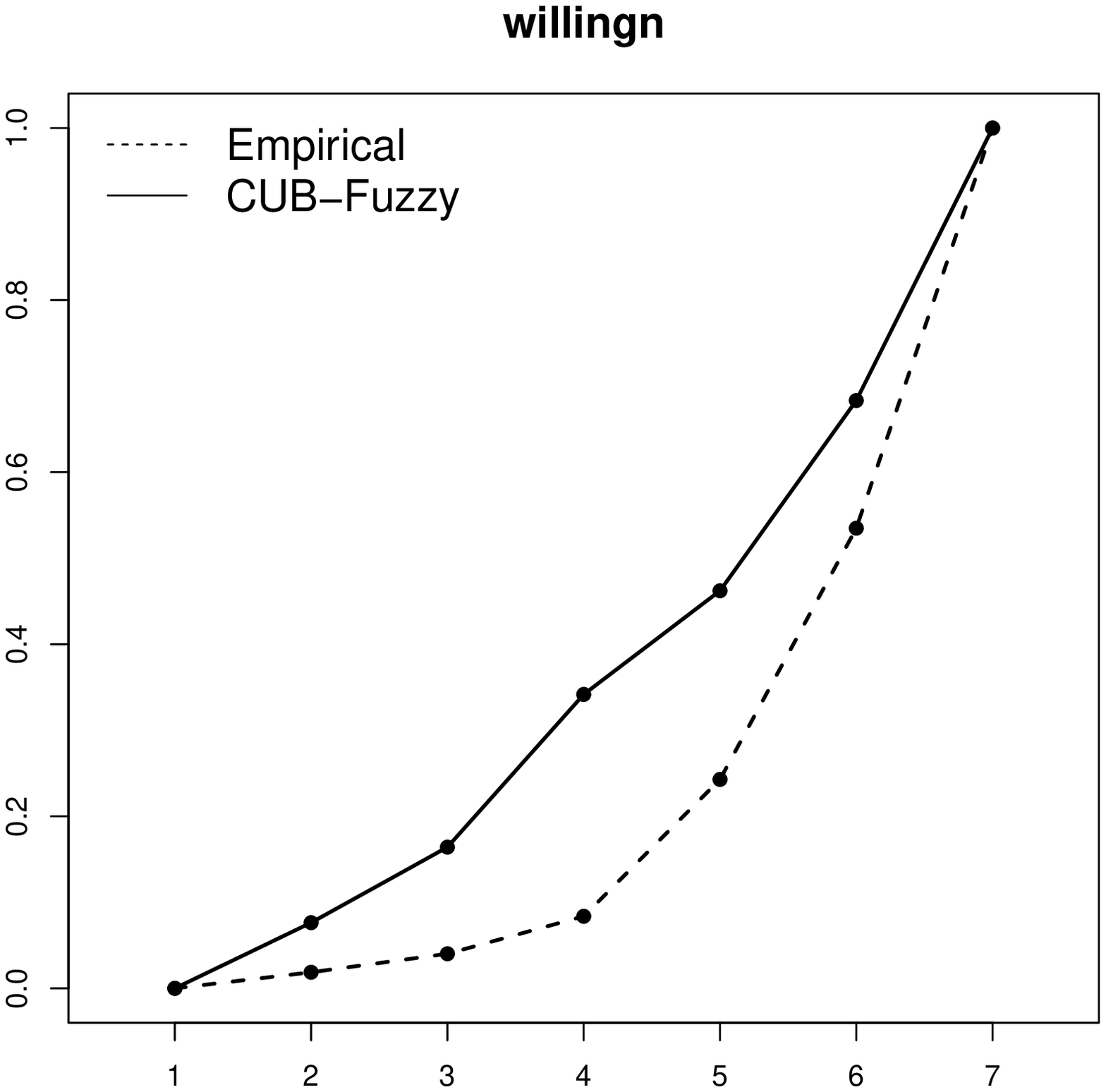}
\includegraphics[scale=0.30]{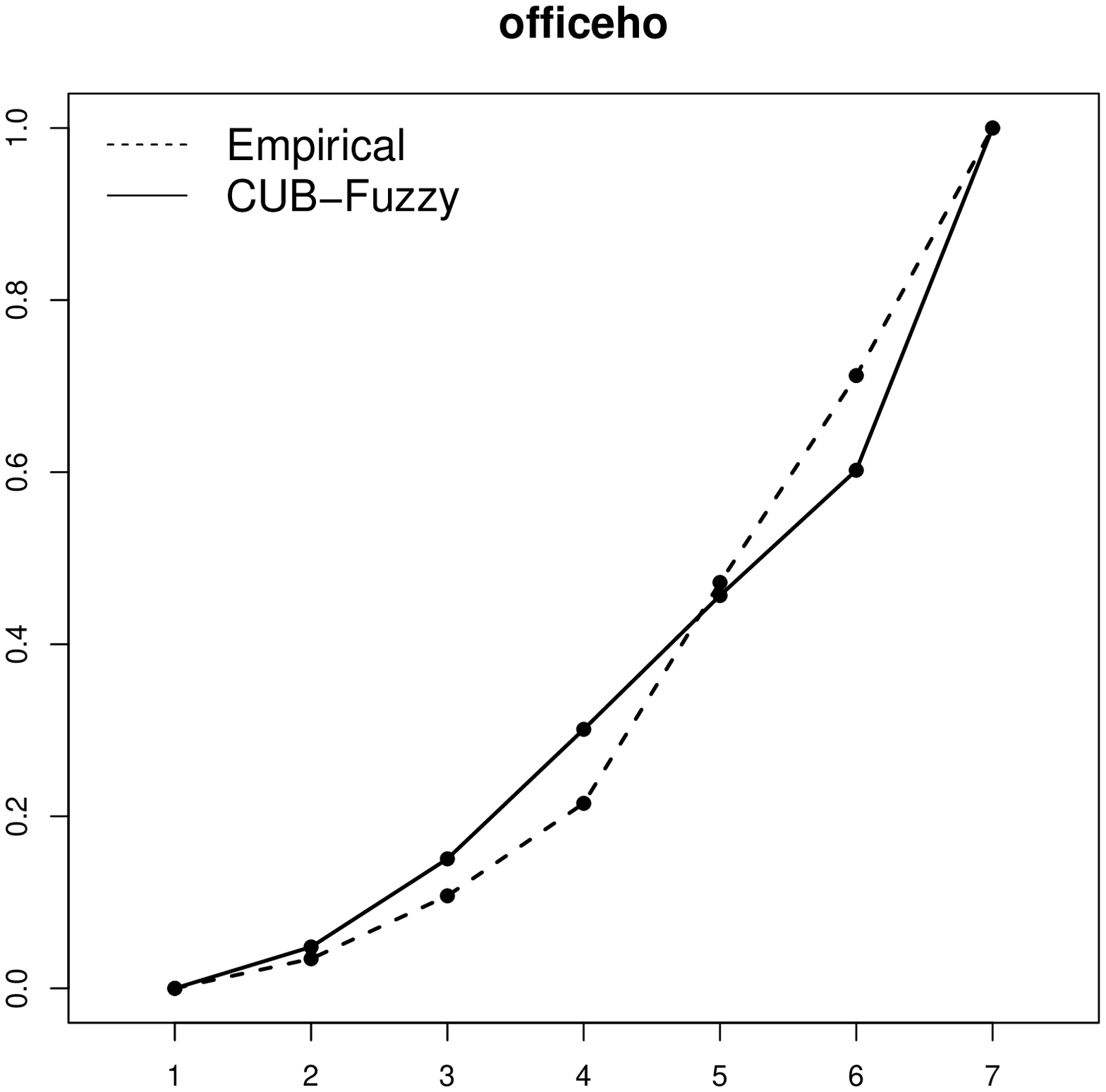}}
\centering{\includegraphics[scale=0.30]{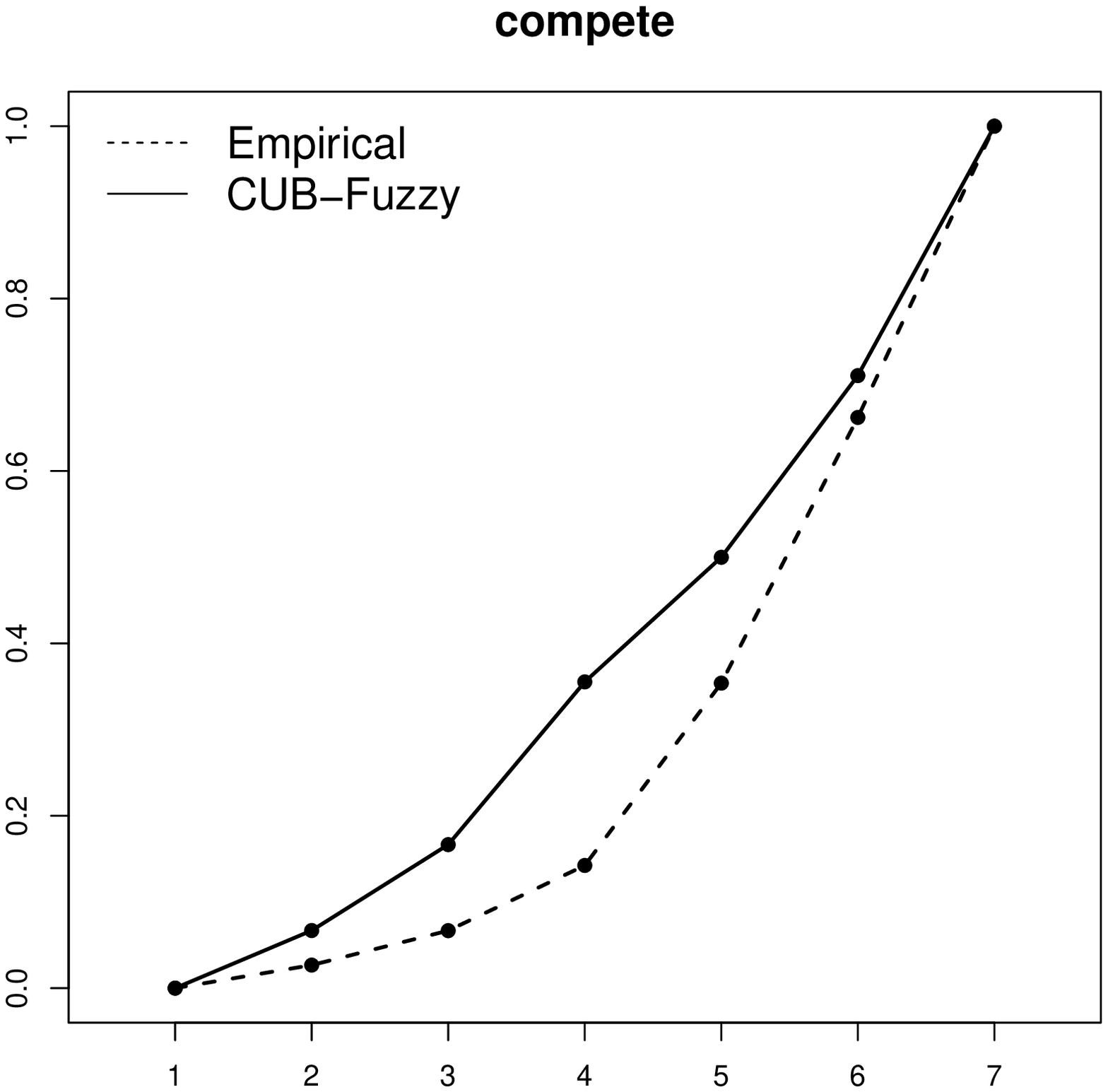}
\includegraphics[scale=0.30]{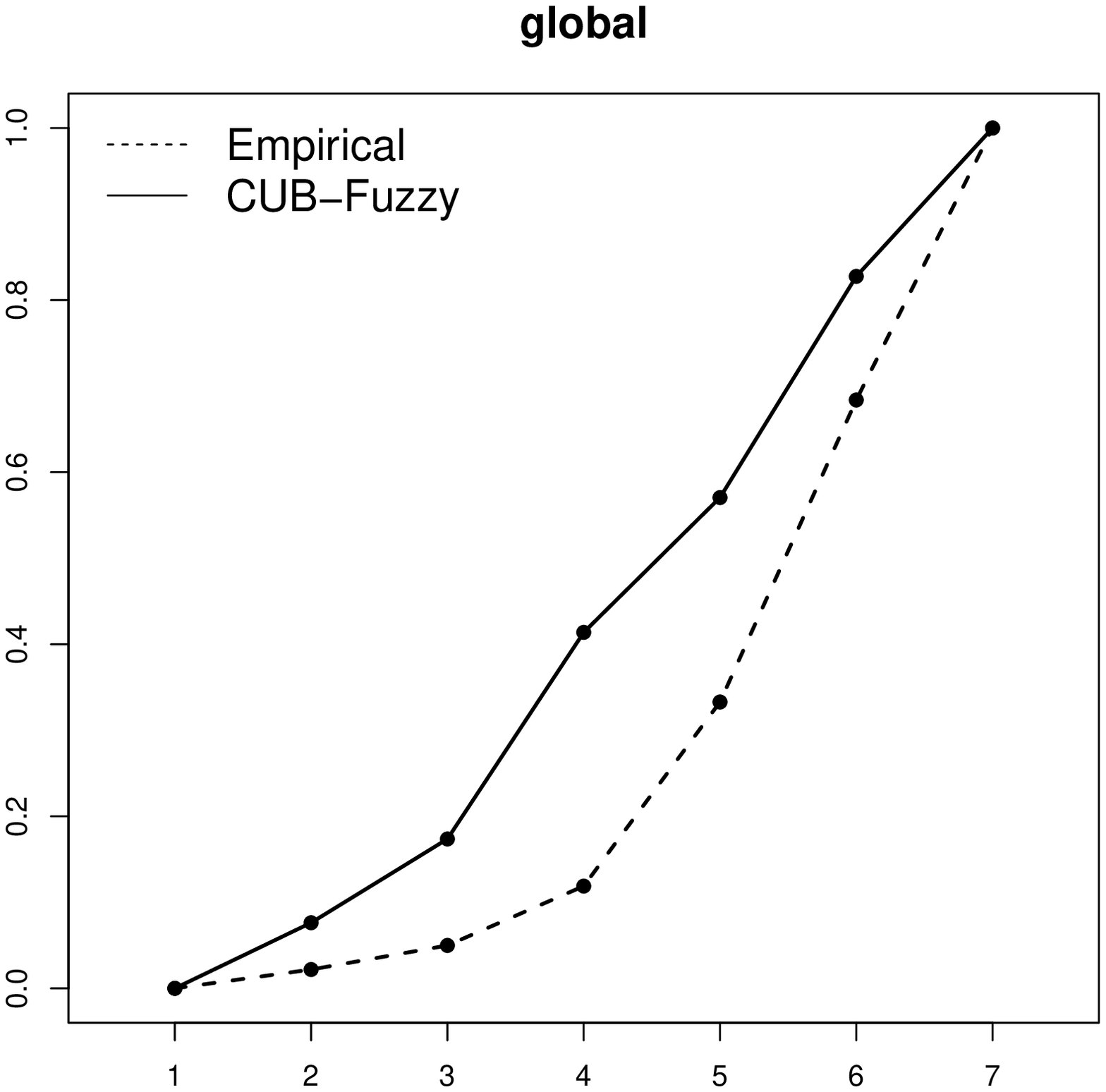}
}
\caption{{\small Comparison of membership functions for the \cub-Fuzzy model (solid line) vs the empirical model (dashed line). }} \label{fig:fig1}
\end{figure}

Figure \ref{fig:fig1} plots the membership function of the \cub-fuzzy evaluation system \eqref{CUBspline} against the empirical one \eqref{Zani} for each item.

For the \cub-Fuzzy approach, the higher the value of $\hat{\pi}_1$ is, the faster the membership degrees increase moving from the indifference point to the maximum of the scale. In addition, notice that the two methods behave quite differently, especially for \texttt{willingn} and \texttt{global} having with the lowest estimated heterogeneity ($ \hat{\pi}_2 = 0.123,$ and $\hat{\pi}_2 = 0.125$, respectively): this indicates a more prominent attitude of the \cub-Fuzzy membership function to correctly discriminate among different levels of heterogeneity, also when moderate.

In order to aggregate membership values, Table \ref{tab:weights} shows the two systems of weights $\{w_k\}$ employable in computing the aggregator index \eqref{(crispZani)}: we refer to \eqref{pesiZani} paired with \eqref{FuzzyPropUF} and with \eqref{FuzzyProp} for the \cub-Fuzzy model (dotted line) and for the empirical model (solid line), respectively. For comparative purposes,  normalized variables loadings derived for the first principal component (PCA1) are also reported. 
\begin{table}[h!]
\caption{Weights systems}   
\label{tab:weights}        
\centering
\begin{tabular}{lccccc}
\hline  & \texttt{informat} & \texttt{willingn} & \texttt{officeho} & \texttt{compete} & \texttt{global} \\
\hline  $w_k$ based on \eqref{FuzzyProp} & 0.216& 0.152 &0.287& 0.196 & 0.146 \\
$w_k$ based on \eqref{FuzzyPropUF} &0.183 &0.212 &0.154 &0.198 &0.254\\ 
\hline 
\hline PCA1  & 0.209   & 0.201  &  0.138  &  0.216   & 0.236  \\
\hline \hline 
\end{tabular}
\end{table}

For the \cub-Fuzzy evaluation system, even if at aggregated level results do not substantially vary at aggregated level  for different weights (see Section \ref{sec:sensitivity}), it turns out that weights based on \eqref{FuzzyProp} do not suitably penalize \texttt{officeho} and \texttt{informat} having a weak importance due to the highest observed uncertainty among the items (for example $1-\hat{\pi_1}=0.398$ for \texttt{officeho} in Table \ref{tab:param_she}). Instead, for the weights based on \eqref{FuzzyPropUF}, the lowest value is attained exactly for \texttt{officeho} ($w_3 = 0.154$). Notice that \texttt{willingn} is assigned a higher weight than \texttt{informat} (w.r.t. the \cub-Fuzzy system of weights), though it shows a higher uncertainty and a lower feeling, comparatively. This is explained by \texttt{willingn} having the most prominent shelter effect at $c=7$, indicating a strong tendency of the distribution to be concentrated at higher categories. Instead, the shelter at $c=5$ for \texttt{informat} acts by deflating the weight of importance (to assess membership, non-membership, etc.) since it is closer to the indifference point and thus tends to penalize a positive evaluation of satisfaction. 

\begin{figure}[h!]
\centering
\includegraphics[scale=0.6]{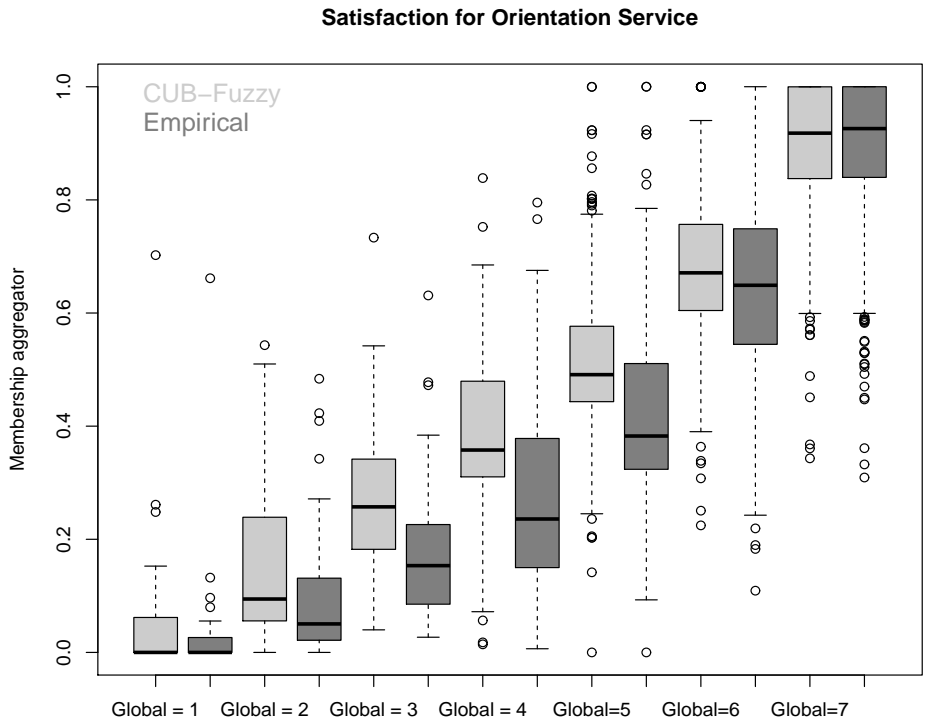}   
\caption{Boxplots of aggregated membership for increasing levels of \texttt{global} satisfaction: comparison between the \cub-Fuzzy (light grey)  and Empirical system (dark grey)}\label{fig:boxplot1}
\end{figure}

Figure \ref{fig:boxplot1} shows as the \cub-Fuzzy system scales membership at aggregated level  more coherently when compared to the expressed \texttt{global} satisfaction. Specifically, we have run both a \cub-Fuzzy evaluation and an empirical system on the first 4 items, omitting \texttt{global} satisfaction. Then, we have stratified the aggregated membership values for the two fuzzy systems across increasing level of \texttt{global} satisfaction (by choosing weights accordingly: the inverse fuzzy proportion of uncertainty for the \cub-Fuzzy proposal and the inverse fuzzy proportion of membership for the empirical one). Then, read on the $y$-axis, it appears that the \cub-fuzzy proposal is more convincing in aggregating the information withdrawn from the first four items if used as a proxy of the global satisfaction. In other words, aggregated values for membership under the \cub-Fuzzy scheme are increasingly more consistent with increasing level of global satisfaction.

\subsection{IFS: CUB-Fuzzy versus spline}\label{sec:CUB_Mar}
\begin{figure}[h!]
\centerline{\includegraphics[width=6.5cm,height=6.5cm]{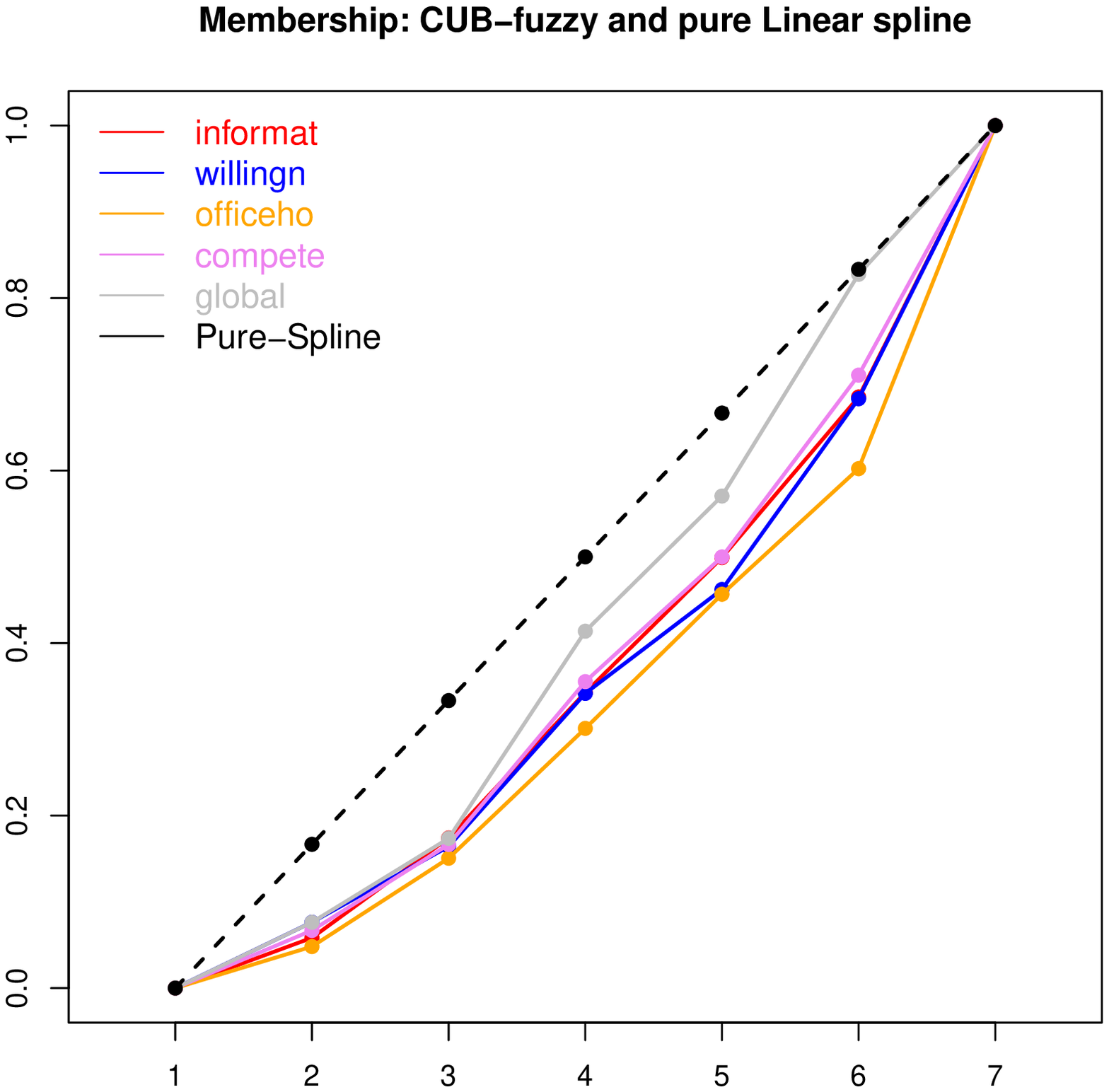}
\includegraphics[width=6.5cm,height=6.5cm]{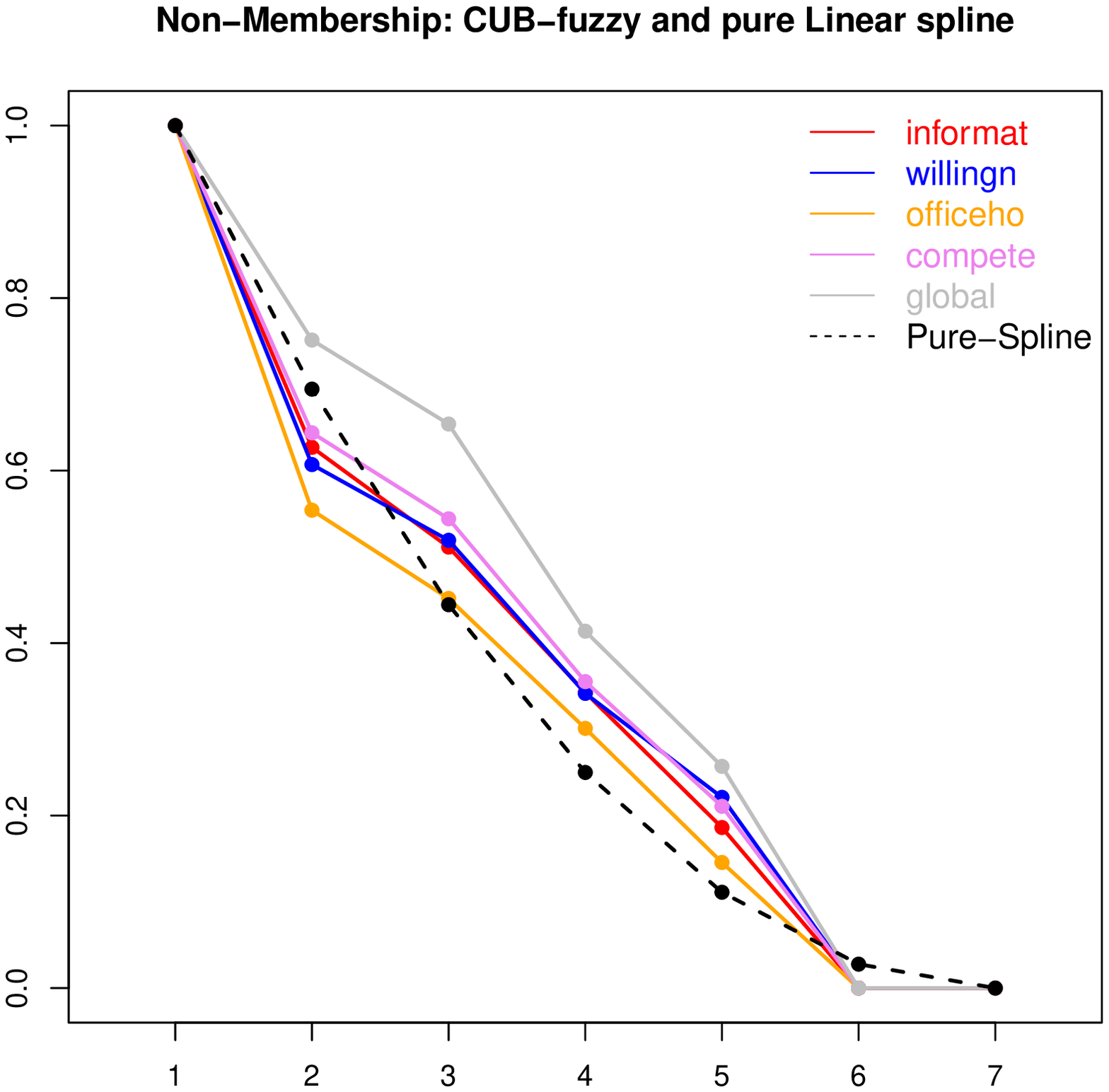}}
\caption{\cub-Fuzzy and spline membership/non-membership functions}\label{figura_nmf}
\end{figure}

Without any prior assessment of experts on specific values for the parameters, membership \eqref{mfspline}, fuzzy uncertainty \eqref{uncMarasini} and non-membership function \eqref{nonmemmar} will be equal for all items, as shown in Table \ref{tab:mf}, where $\epsilon = 1,\, a=1, b=m-1$ and $\theta = \eta$ have been set to account for the balanced scale as recommended by \cite{Marasini}. In particular, in the following we have tested the \cub-Fuzzy proposal against the spline uncertainty with $\theta = \eta = 1$ due to its interpretation as a risk measure: indeed, in this case the spline uncertainty fuzzy function corresponds to the variance of a Bernoulli random variable whose success trial (the membership to $A$) has probability of occurrence set to $\mu_A(r)$. Dually, the fuzzy uncertainty prescribed by the \cub-Fuzzy system accounts for risk in terms of heterogeneity.
\begin{table}
\caption{Spline Fuzzy functions}
\label{tab:mf}        
\centering
\scalebox{1}{\begin{tabular}{lccccccc}
\hline \hline  &  $R = 1 $ & $R = 2$ & $R = 3$ & $R = 4$ & $R = 5$ & $R = 6$ & $R =7 $ \\ \hline \hline
Membership degree     & 0  & 0.10  & 0.30     & 0.50     & 0.70   & 0.90  &  1  \\ 
Non-membership degree & 1  & 0.81 & 0.49    & 0.25    & 0.09  & 0.01 &  0  \\ 
Uncertainty degree  & 0  & 0.09 & 0.21    & 0.25    & 0.21  & 0.09 &  0  \\
\hline \hline  
\end{tabular}}
\end{table}
In Figure \ref{figura_nmf}, such spline membership and non-membership values are compared with those obtained with the \cub-Fuzzy proposal. Observe that the membership given in Definition \ref{CUBspline} is more accurate and naturally shaped to the data. 
What is constant over categories in the \cub-Fuzzy model is the Fuzzy uncertainty function in Definition \ref{unc} (that is $1-\hat{\pi}_1$ in Table 
\ref{tab:param_she}). Indeed $1-\hat{\pi}_1$ (equivalently, $1-\hat{\pi}$ when \textit{shelter effect} is not significant) measures the overall \textit{fuzziness}, independently from the membership and non-membership degrees; dually, $\hat{\pi}_1$ quantifies the level of accuracy of the Fuzzy evaluation system. This feature is only partially accomplished by the spline method as the Fuzzy uncertainty function has symmetric values around the indifference point, see Table \ref{tab:mf}. Notice that in general, the uncertainty function of an IF evaluation system is set in such a way that it attains a maximum at the indifference category if available. Instead, we consider the uncertainty as uniformly spread along the scale, so that it can be considered as a feature of the item and not of a single category.

Comparisons between the two methods based, for example, on the composite indicator \eqref{(crispZani)} are meaningless as for the spline one since $\tilde{\mu}_A(r) = \mu_A(r) = 
\mu_A^{(k)}(r)$ for all $k$. For this reason, we propose to aggregate the $k$-th item uniformly across respondents, achieving a complete IFS evaluation system. More specifically, by keeping the notation introduced in Section \ref{sec:uncertainty}, for the membership function we compute
$\bar{\mu}^{(k)} = \frac{1}{n}\sum\limits_{j=1}^n \mu_A^{(k)}(r_{j,k})$ both for the (linear) spline and the \cub-Fuzzy evaluation systems, see Table \ref{table:compareitem}. The same is done for the non-membership functions as well as the Fuzzy score and accuracy measures. 
\begin{table}[h!]
\caption{Fuzzy functions aggregated per item}\label{table:compareitem}
\centering
\scalebox{0.9}{
\begin{tabular}[pos=l]{llccccc}
\hline \hline Function & Method $\diagdown$ Item: & \texttt{informat} & \texttt{willingn} & \texttt{officeho} & \texttt{compete} &  \texttt{global} \\
\hline  \multirow{2}{*}{Membership}  & \cub-Fuzzy & 0.651& 0.743 &0.571& 0.681& 0.752\\ 
& Spline & 0.780& 0.853 &0.733 &0.797 &0.814\\   
\hline  \multirow{2}{*}{Non-membership}  &  \cub-Fuzzy &0.136& 0.091& 0.157 &0.131 &0.132 \\
& Spline &  0.113 &0.072 &0.151 &0.104& 0.086\\
\hline  \multirow{2}{*}{Fuzzy uncertainty}  &  \cub-Fuzzy &0.212 &0.166& 0.272& 0.188 &0.116 \\
& Spline &0.107& 0.075& 0.116 &0.099 &0.100\\
\hline  \multirow{2}{*}{Fuzzy score}  & \cub-Fuzzy & 0.515& 0.652& 0.414& 0.550& 0.620\\
& Spline &0.667 &0.781& 0.581 &0.693& 0.728 \\
\hline  \multirow{2}{*}{Fuzzy accuracy}  & \cub-Fuzzy & 0.788 &0.834& 0.728 &0.812& 0.884 \\ 
& Spline &  0.893 &0.925 &0.884 &0.901 &0.900 \\ 
\hline  \hline \\
\end{tabular}}
\end{table}
As we see from Table \ref{table:compareitem}, the spline approach does not sufficiently discriminate the different levels of uncertainty among the items, yielding to a narrow range for both the fuzzy accuracy and uncertainty. Conversely, the \cub-Fuzzy proposal offers a major flexibility in grading Fuzzy indicators according to the observed uncertainty. 
In particular, we stress that the accuracy of the \cub-Fuzzy approach is penalized for items with higher global uncertainty in the sense of \cub models, while it increases for items corresponding to a weaker indeterminacy. For instance, with reference to Table \ref{tab:param_she}, the maximum estimated overall indeterminacy corresponds to \texttt{officeho}  ($1-\hat{\pi}_1 = 0.398$), whereas the minimum corresponds to \texttt{global} ($1-\hat{\pi}_1 = 0.172$). As a result, these items coherently are assigned the minimum and maximum levels of the accuracy, respectively, while the spline model is more restrictive in accounting for this variability.

\begin{figure}[h!]
\centering

\includegraphics[scale=0.55]{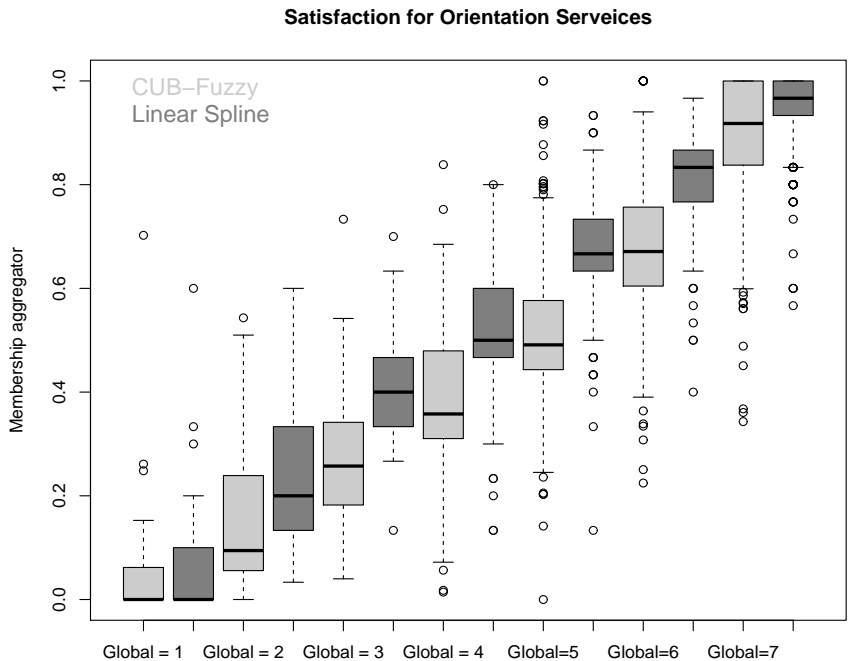} %
\caption{Boxplots of aggregated membership for increasing levels of \texttt{global} satisfaction: comparison between the \cub-Fuzzy (light grey) and linear spline systems (dark grey) }\label{fig:boxplot2}
\end{figure}

\begin{figure}[h!]
\centering
\includegraphics[scale=0.6]{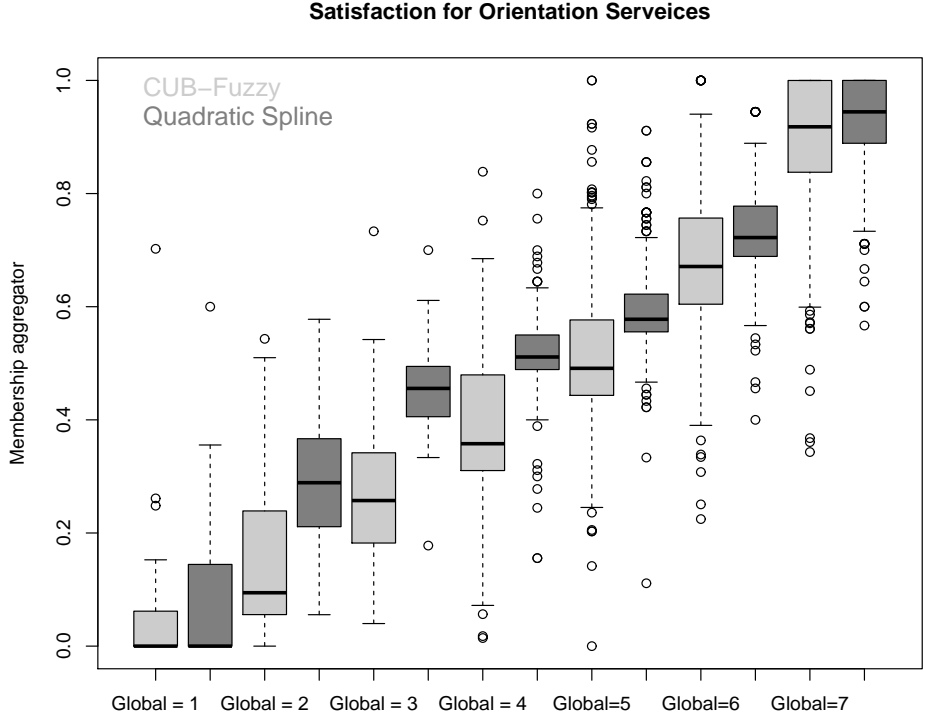}
\caption{Boxplots of aggregated membership for increasing levels of \texttt{global} satisfaction: comparison between the \cub-Fuzzy (light grey) and quadratic spline systems (dark grey)}\label{fig:boxplot3}
\end{figure}

Figures \ref{fig:boxplot2}-\ref{fig:boxplot3} show as the \cub-Fuzzy system scales -on  the $y$ axis- membership at aggregated level in a comparable way as the linear spline  conditional on increasing scores for \texttt{global} satisfaction, whereas it is more adequate than the quadratic spline system. Spline memberships have been aggregated with uniform weights across items. Comparable results are obtained if normalized variables loadings derived from factor analysis are considered.\\

Finally, the ultimate step of the Fuzzy evaluation procedure is to provide a measure of the overall uncertainty by aggregating the functions 
given in Table \ref{table:compareitem} among items: as example, for the membership function we compute $\bar{\mu} = \sum_{k=1}^K w_k \bar{\mu}^{(k)}$
with $\bar{\mu}^{(k)}$ given in Table \ref{table:compareitem} (first row); similarly, for all fuzzy functions. For the \cub-Fuzzy method we employ the weights \eqref{pesiZani} paired with \eqref{FuzzyPropUF}. For the spline method, instead, we shall consider uniform weights also across the items \cite{Marasini}. 
The resulting Fuzzy composite aggregators are reported in Table \ref{tab:compare}. 
\begin{table}[h!]
\centering
\caption{Fuzzy composite indicators (aggregation of items)}\label{tab:compare}
\scalebox{0.8}{
\begin{tabular}{llccccc}
\hline \hline & Weights & Membership & Non-membership & Uncertainty & Score & Accuracy \\
\hline  \cub - Fuzzy & \footnotesize{Log inverse of Fuzzy unc.} & 0.690 & 0.128 &0.182& 0.562 &0.818 \\ 
\hline  Spline &  \footnotesize{Uniform weights} & 0.795 & 0.105 & 0.099 & 0.690& 0.900\\
\hline  \\                                                                            
\end{tabular}}
\end{table}
Thus, we can conclude that the proposed evaluation system based on \cub models is safer in assigning fuzzy values, being designed to account for heterogeneity and stylistic responses in the data. Nevertheless, it does not miss to provide a global positive picture (in terms of membership and accuracy).

In conclusion, in order to disclose the different perspectives offered by the \cub-Fuzzy analysis of questionnaire, we combined a two-component PCA analysis with a $k$-means clustering on the data-matrix according to a so-called \textit{tandem scheme} \cite{ArabieHubert}. Specifically, we compared the closeness of the derived classification with that of a $k$-means algorithm applied to the IWAM aggregators \eqref{IWAM}:  $k=5$ was set for the $k$-means algorithm to identify certainly unsatisfied ($R=1$), fairly unsatisfied ($R=2,3$), uncertain ($R=4$), fairly satisfied ($R=5,6$), and certainly satisfied ($R=7$) respondents.  Table \ref{tab:cohen} reports the Cohen's $\kappa$ measure to assess agreement between the two corresponding classifications, along with lower and upper confidence bounds:
\begin{table}[h!]
\centering
 \caption{Cohen's $\kappa$: agreement between the classification obtained from $k$-means on the first two PCA components and that obtained from $k$-means on membership and non-membership aggregators \eqref{IWAM} for different fuzzy methods. }\label{tab:cohen}
\begin{tabular}{lccc}
\hline & Lower & Estimate & Upper \\
\hline  \cub-Fuzzy & 0.45 & 0.49 & 0.53 \\
\hline Linear Spline & 0.33 & 0.36 & 0.38 \\
\hline Quadratic Spline & 0.18 & 0.21 & 0.24 \\
\hline 
\end{tabular}
\end{table}

Results do not substantially vary if considering the polychoric correlation to run PCA. Notice that fuzzy clustering has a precise meaning in the literature \cite{Everitt}, which is not involved in the present analysis: this perspective will be the subject of future investigation. 

\subsection{Sensitivity Analysis}\label{sec:sensitivity}
A sketch of sensitivity analysis is here accomplished to validate the \cub-Fuzzy proposal against the other fuzzy alternatives here considered.

\begin{figure}[h!]
\centering
\includegraphics[scale=0.55]{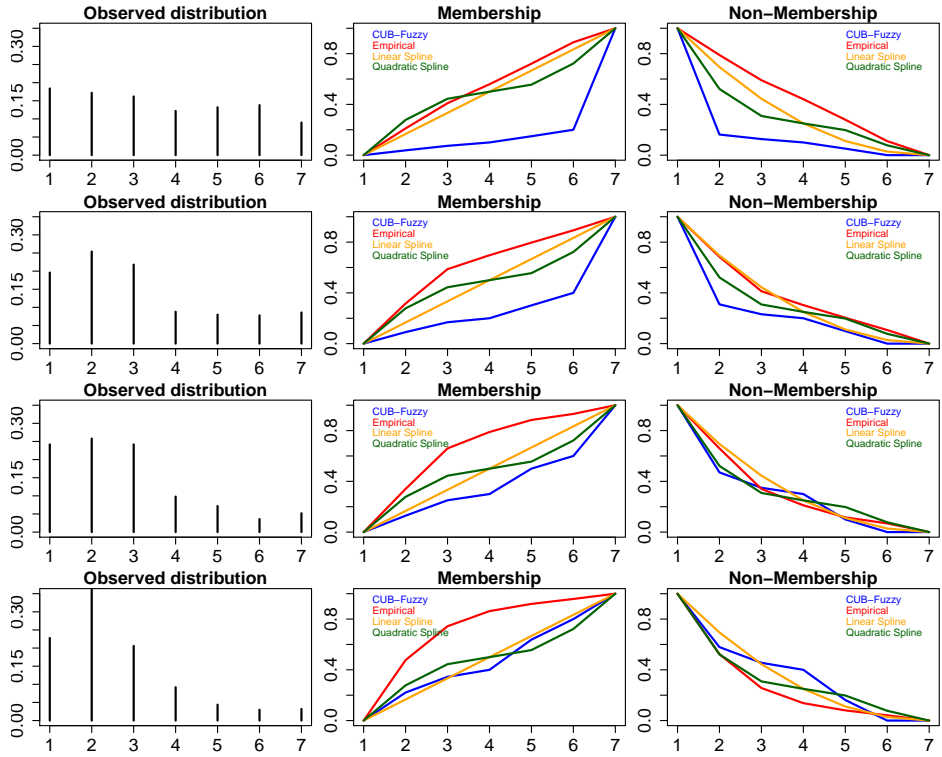}
\caption{Comparison of Membership and Non-Membership for right-tailed rating distributions ($\xi = 0.8$) corresponding to decreasing levels of heterogeneity: $\pi = 0.2$ (top) up to $\pi = 0.8$ (bottom)}\label{fig:comp1}
\end{figure}

\begin{figure}[h!]
\centering
\includegraphics[scale=0.55]{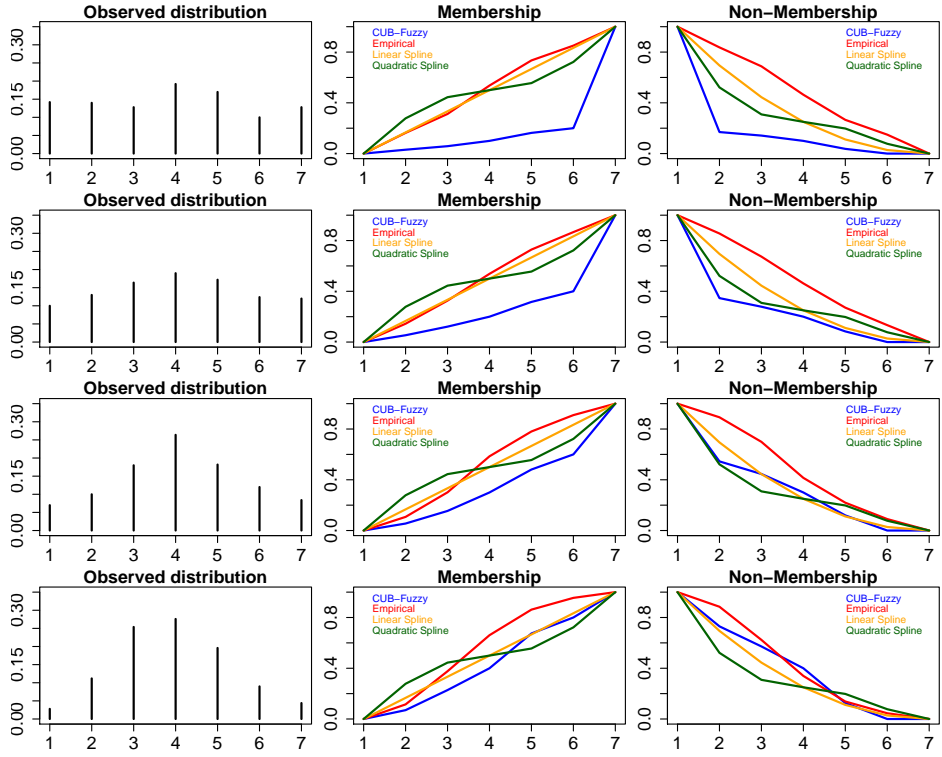}
\caption{Comparison of Membership and Non-Membership for symmetric rating distributions ($\xi = 0.5$ for decreasing levels of heterogeneity: $\pi = 0.2$ (top) up to $\pi = 0.8$ (bottom)}\label{fig:comp2}
\end{figure}

\begin{figure}[h!]
\centering
\includegraphics[scale=0.55]{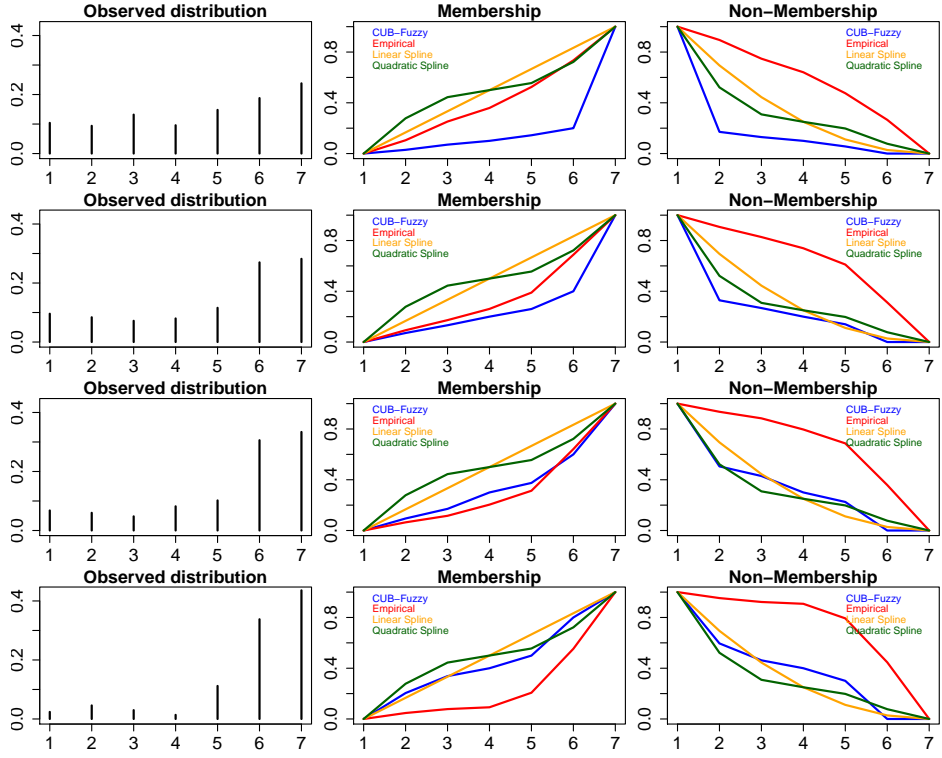}
\caption{Comparison of Membership and Non-Membership for left-tailed rating distributions ($\xi = 0.1$) for decreasing levels of heterogeneity: $\pi = 0.2$ (top) up to $\pi = 0.8$ (bottom)}\label{fig:comp3}
\end{figure}

First, read top to bottom, Figures \ref{fig:comp1}-\ref{fig:comp3} display membership and non-membership functions for decreasing levels of heterogeneity and for right-tailed, symmetric and left-tailed rating distributions generated from varying \cub distributions: $\xi = 0.8, \xi = 0.5, \xi = 0.1$, respectively, and for each of them $\pi= 0.2$ (top), $\pi = 0.4, \pi=0.6, \pi=0.8$ (bottom). Imagine that the distributions correspond to measurements on a scale $1$=``\textit{extremely dissatisfied}'' up to $7$=``\textit{completely satisfied}''. Then Figures \ref{fig:comp1}-\ref{fig:comp3} correspond to overall dissatisfied, overall neutral and overall satisfied respondents, respectively. It appears evident that the \cub-Fuzzy proposal, being naturally shaped to the data, is safer and more integrated with the observed scores. In addition, for data with low heterogeneity (the bottom panels), it is globally intermediate between the linear and quadratic splines. The price of estimation - which can be promptly run by means of the \texttt{R} package `CUB' \cite{CUBpackage} - is compensated with no need of prior choice for  spline degrees and parameters, and  with a more flexible and versatile tool for uninformative circumstances.


Secondly, we assess distances between aggregated Intuitionistic Fuzzy sets obtained with different weighting and fuzzy systems. Specifically, the (normalized) Hamming distance between IFS sets   is computed \cite{FuzzyDistance}. Briefly, if 
$$ B = \{<x, \mu_B(x),\nu_B(x)>| x \in X \}, \qquad C = \{<x, \mu_C(x),\nu_C(x)>| x \in X \}$$
 are two IF evaluation systems defined for the universe of discourse $X$, then the (normalized) Hamming distance between $B$ and $C$ is defined as:
$$ d_H(B,C) = \dfrac{1}{2n}\sum_{i=1}^n \bigg(|\mu_B(x_i) - \mu_C(x_i)| + |\nu_B(x_i) - \nu_C(x_i)| + |u_B(x_i) - u_C(x_i)|\bigg),$$
with $n$ being the number of observations. 
Back to our case study, in order to support the consistency of the \cub-Fuzzy proposal, we will compute the distance between the aggregated IWAM fuzzy sets \eqref{IWAM} corresponding to different choices of the weighting system.
For $l=1,2$, let $\bm w^{(l)} = \{w_1^{(l)},\dots,w_K^{(l)}\}$ be two alternative choices for weights for which:
$$ \tilde{\mu}_A^{(l)}(\bm r_i) = \sum_{k=1}^K w_k^{(l)}\,\mu_A^{(k)}(r_{i,k}), \qquad \tilde{\nu}_A^{(l)}(\bm r_i) = \sum_{k=1}^K w_k^{(l)}\,\nu_A^{(k)}(r_{i,k})$$
are the corresponding IWAM for membership and non-membership functions across $K$ items. Then, let:
$$ B = \{<\bm r_i, \tilde{\mu}_A^{(1)}(\bm r_i),\tilde{\nu}_A^{(l)}(\bm r_i) > | i=1,\dots,n\}$$ 
$$ C = \{<\bm r_i, \tilde{\mu}_A^{(2)}(\bm r_i),\tilde{\nu}_A^{(2)}(\bm r_i) > | i=1,\dots,n\}.$$

\begin{table}[h!]
\caption{Normalized Hamming distance between IWAM aggregators \eqref{IWAM}  of the \cub-Fuzzy system with different weights.}\label{tab:dist1}
\centering
\begin{tabular}{l|ccc}
\hline $d_H(B,C)$ & \multicolumn{3}{c}{$\bm w^{(2)}$}\\
\hline  $\bm w^{(1)}$ & \eqref{pesiZani} with \eqref{FuzzyProp} & PCA loadings & Uniform \\ 
\hline \eqref{pesiZani} with \eqref{FuzzyPropUF} & 0.028 & 0.006 & 0.011 \\
\hline 
\end{tabular}
\end{table}

As reported in Table \ref{tab:dist1}, the restrained values of the distance support the relative indifference between the fuzzy proportion of uncertainty and membership values according to \eqref{pesiZani} for the \cub - Fuzzy Proposal, and with weights derived from PCA-type procedures. Nevertheless, we acknowledge that weights based on the fuzzy proportion of uncertainty are the most natural and always applicable choice for \cub-Fuzzy systems, and that weights based on PCA-type are an acceptable option only if the first component explains an appreciable amount of variability. In addition, some of the \cub-Fuzzy indicators are stable with respect to the reversion of the scale, that is for samples $m - r_{j,k} +1$ with $j=1,\dots,n$ and $k=1,\dots,K$. Indeed, the \cub random variable $R$ is reversible: if $R \sim $ \cub($\pi$, $\xi$) over the $m$ categories, then $m - R +1  \sim $ \cub($\pi$, $1-\xi$) \cite{DelPic2005}. If the scale is balanced and $R \sim$ \cub($\pi, \xi)$, with shelter at $c$ measured by $\delta$, then $m-R +1 \sim$ \cub($\pi, 1-\xi$) with shelter at $m-c+1$ measured by $\delta$ as well. Due to this property, for reversed ratings the \cub-Fuzzy accuracy and uncertainty functions remain unchanged, differently from the membership and non-membership degrees, being dependent on the distribution functions.  Instead, the spline Fuzzy functions would change only at an aggregated level. In addition, also weights \eqref{pesiZani} paired with \eqref{FuzzyPropUF} are invariant, being directly related to \cub uncertainty parameters. Conversely, if paired with the Fuzzy proportion of membership degrees \eqref{FuzzyProp}, such weights would vary. f

\section{Comments and conclusions}
The present contribution fosters the application of \cub mixture models for ordinal rating data to account for uncertainty of choices within a fuzzy analysis of questionnaires. The resulting \cub-Fuzzy procedure is well-suited for broad applications since it allows to deal with veracity of rating by means of the assessment of membership and non-membership grades.

From the statistical modelling point of view, the proposal sheds new light on the \cub paradigm, conveying its vague definition of uncertainty into a precise frame. From the point of view of fuzzy analysis, it roots membership and non-membership assessments on the basis of sound statistical procedures, thus fuzzy functions gain reliability. Moreover, the proposal is build in such a way that the data structures and relations are preserved, and -when pertinent- results match with traditional methods.

The proposal stems from the spline approach introduced in \cite{Marasini}, suitably adjusted with the \cub uncertainty parameter to let the spline approach be more insightful and freed of subjectivity of parameters choice. The methodology here introduced meets also a classical proposal based on the empirical distribution function \cite{Cheli,Zani}, which is adjusted within a general IFS framework. The procedure provides a refined tool to account for heterogeneity and other forms of nuisances, as meant by the rationale of \cub uncertainty.  In particular, the \cub model uncertainty $1-\pi$ is equal to the Fuzzy hesitancy level for each item, and the accuracy function results to be more sensitive to different sources of indeterminacy as heterogeneity and \textit{shelter effect} (see Table \ref{tab:compare}). As a result, the \cub uncertainty measure is validated as an effective Fuzzy composite indicator.  

Summarizing, the spline methods for fuzzy analysis of questionnaire is a valuable methodology, whose main criticisms are the subjectivity of choices and the lack of a statistical foundation, where diversity in scale-usages needs to be taken into account. The \cub-Fuzzy system overcomes these pitfalls by defining fuzzy functions anchored to the empirical distribution functions and adjusted for uncertainty in the data. Encoding the uncertainty degrees implies that the resulting \cub-Fuzzy system is freed of subjectivity of parameters values for spline functions, without the need of choosing between a linear or a quadratic splines since the vagueness of responses induced by the scale is automatically charged by the uncertainty parameter. In addition, it is grounded on ML estimation and in this sense it is more robust. In the same vein, in order to let the \cub-fuzzy analysis of questionnaire adhere to respondents' subjectivity and not to that of scholars and judges,  the uncertainty parameter can be estimated on subjective basis ($\pi_i$) and linked to responses drivers (covariates $Y_i$) by means of a logistic transform:
$$ logit(\pi_i) = \beta_0 + \bm \beta^{'} Y_i. $$ 
Future developments in this direction will be adressed to take into account the Hesitant Fuzzy Set framework too, as proposed in \cite{Marasini,Torra}. Nevertheless, it is worth to underline that the \cub-fuzzy proposal can be enriched by the specification of covariates to disclose response profiles, but it is valid \textit{per se}, conversely to some other traditional methods.

In conclusion, let us underline that the proposal does not advance a brand new statistical model, rather it is tailored to boost the idea that any statistical model should be prone to offer some veracity analytics of data.

Traditional models are able to discriminate sharply between satisfied and dissatisfied respondents, at the price of more involved model-specifications and no measure of uncertainty (and, dually, no measure of accuracy/reliability of ratings).  Conversely, the proposal offers a multifaceted interpretation of results, with both \textit{local} modelling (as for cumulative and partial credit models, for instance), and  global assessments (feature that is inherited by \cub models). For instance, it allows discrimination of items in terms of their capabilities of identifying satisfied and unsatisfied respondents (and, in full generality, members and not-members of latent classes). This advantage can be appreciated the more data are heterogeneous.

\end{document}